\documentclass[%
]{mpi2015-cscpreprint}

\usepackage{epstopdf}
\usepackage[caption=false]{subfig}
\usepackage{amsmath,graphicx}

\newcommand{\diagdots}[3][-25]{%
	\rotatebox{#1}{\makebox[0pt]{\makebox[#2]{\xleaders\hbox{$\cdot$\hskip#3}\hfill\kern0pt}}}%
}
\usepackage{amsfonts}
\usepackage{amssymb}
\usepackage[english]{babel}
\usepackage{pgfplots}
\usepackage{nomencl}
\makenomenclature
\usetikzlibrary{matrix}
\usepackage{parskip}
\usepackage{mymacros}
\usepackage{amsthm}
\usepackage{tikz}
\usepackage{enumerate}
\usetikzlibrary{calc,patterns,decorations.pathmorphing,decorations.markings}

\usepackage{cite}
\usepackage{esvect} 
\usepackage{soul} 
\usepackage{tikz}
\tikzstyle{vertex}=[circle, draw, inner sep=0pt, minimum size=6pt]
\usepackage{graphicx}
\usepackage{float}
\usepackage{mathtools}
\usepackage{todonotes}
\usepackage{booktabs}
\usepackage{ifdraft}
\usepackage{algorithm}
\usepackage{algpseudocode}
\usepackage{setspace}
\usepackage{color}

\usepackage[normalem]{ulem}

\renewcommand{\hat}{\widehat}
\renewcommand{\tilde}{\widetilde}

\usepackage[square,sort&compress,comma,numbers]{natbib}
\bibpunct[, ]{[}{]}{,}{n}{,}{,}
\newcounter{mymac@matlab}
\newcommand{\matlab}{MATLAB%
	\ifnum\value{mymac@matlab}<1%
	\textsuperscript{\textregistered}%
	\setcounter{mymac@matlab}{1}%
	\fi%
}

\theoremstyle{plain}

\newtheorem{remark}{Remark}

\theoremstyle{definition}

\theoremstyle{remark}

\usepackage{graphicx}

\usepackage{amssymb}
\usepackage{amsthm}

\begin{document}
  

\title{Semi‐active damping optimization of vibrational systems using the reduced basis method}
  
\author[$\ast~1$]{Jennifer Przybilla}
 \affil[$\ast$]{Max Planck Institute for Dynamics of Complex Technical Systems, Sandtorstra{\ss}e 1, 39106 Magdeburg, Germany.}
 \affil[1]{\email{przybilla@mpi-magdeburg.mpg.de}, \orcid{0000-0002-8703-8735}}
 \affil[2]{\email{pontes@mpi-magdeburg.mpg.de}, \orcid{0000-0001-6433-6142}}
 \affil[3]{\email{benner@mpi-magdeburg.mpg.de}, \orcid{0000-0003-3362-4103}}
  	
\author[$\,\ast~2$]{Igor Pontes Duff}  
\author[$\ast~3$]{Peter Benner}
  
\shorttitle{Damping optimization of vibrational systems}
\shortauthor{J. Przybilla, I. Pontes Duff, P. Benner}
\shortdate{}

\title{Semi-active damping optimization of vibrational systems using the reduced basis method}

\abstract{%
In this article, we consider vibrational systems with semi-active damping that are described by a second-order model.
In order to minimize the influence of external inputs to the system response, we are optimizing some damping values.
As minimization criterion, we evaluate the energy response, that is the $\cH_2$-norm of the corresponding transfer function of the system.
Computing the energy response includes solving Lyapunov equations for different damping parameters.
Hence, the minimization process leads to high computational costs if the system is of large dimension.

We present two techniques that reduce the optimization problem by applying the reduced basis method to the corresponding parametric Lyapunov equations. In the first method, we determine a reduced solution space on which the Lyapunov equations and hence the resulting energy response values are computed approximately in a reasonable time. The second method includes the reduced basis method in the minimization process. To evaluate the quality of the approximations, we introduce error estimators that evaluate the error in the controllability Gramians and the energy response.
Finally, we illustrate the advantages of our methods by applying them to two different examples. }

\novelty{
In this paper, we propose methods to reduce the problem of minimizing the energy response of a vibrational system by optimizing corresponding damping gains.
Additionally, we propose error estimators that are mandatory for the methods to evaluate the quality of the approximated controllability Gramian and the corresponding energy response.
}
\keywords{reduced basis method, damping optimization, vibrational systems, error estimation, model order reduction}

\maketitle


\section{Introduction}\label{sec:Intro}
When constructing large civil engineering infrastructure such as buildings or bridges, external vibrational forces like wind perturbations or earthquakes need to be taken into account.
Due to the continuous improvement in engineering construction, which provides for lighter and finer structures, corresponding infrastructures have become more susceptible to large deflections and fatigue when external forces with dominant frequencies close to the natural frequencies of the structure are applied. 
To prevent this effect, dampers are designed in order to remove critical energies from the physical system.
We consider the vibrational system
\begin{align}\label{eq:SOsys}
\begin{split}
M\ddot{x}(t) + D(g)\dot{x}(t) + Kx(t) &= Bu(t),\\
y(t) &= C x(t)
\end{split}
\end{align}
where $M\in\Rnn$ is the mass matrix, $D(g)\in\Rnn$ is the damping matrix and $K\in\Rnn$ the stiffness matrix.
We assume, that $M$ and $K$ are symmetric and positive definite, and that $D(g)$ is symmetric and positive semidefinite for all parameters $g\in\cD$, where $\cD\subseteq\R^{\ell}$ is the parameter set.
Additionally, the matrix $B\in\Rnm$ is the input matrix and $C\in\Rpn$ the output matrix.
The vectors $u(t)\in\Rm$, $x(t)\in\Rn$ and $y(t)\in\Rp$ describe the input, the state and the output of the system, respectively.
The damping of the system consists of two parts: an internal damping of small magnitude and external dampers designed to limit the influence of the input to the output of the system.
Hence, the parameter-dependent damping matrix $D(g)$ is composed of 
$
D(g) = D_{\intern} + D_{\extern}(g),
$
where $D_{\intern}$ describes the internal damping and $D_{\extern}(g)$ describes the external damping.

The internal damping can be modeled in different ways.
In this work we use a multiple of the critical damping  
\begin{align}\label{eq:intDamp}
D_{\intern} := 2\alpha M^{\frac{1}{2}}\left(M^{-\frac{1}{2}}KM^{-\frac{1}{2}}\right)^{\frac{1}{2}}M^{\frac{1}{2}},
\end{align}
for $\alpha\ll 1$, which is a widely used convention, and was applied, for example, in \cite{morBenTT11, morBenTT13,Ves11,morTruV09,PazK18}.
However, our theory holds for every modal damper, i.e., every damping matrix that is simultaneously with $M$ and $K$ diagonalizable.
Other methods to model the internal dampings are presented in \cite{JakMTetal21,KuzTT12}.

External dampers typically depend on two variables, the position and the damping values.
In this paper, we will focus on the optimization of the damping values to attenuate the effects from the external forces. 
Hence, the external damping matrix $D_{\extern}(g)$ depends on the damping gains $g = \begin{bmatrix}
g_1,\dots , g_{\ell}
\end{bmatrix}$, that represent the friction coefficients, and the matrix $F\in\R^{n\times\ell}$ that describes the position of the dampers such that 
\begin{align*}
D_{\extern}(g):=FG(g)F^\T,\qquad G(g) := \diag{g_1, \dots ,g_{\ell}}.
\end{align*}
We assume that the number of dampers $\ell$ is significantly smaller than the dimension $n$.
Additionally, we assume that the damping gains are fixed over time and lie in given intervals $[g_i^-, \, g_i^+]\subset\R_+$, for all $i = 1,\dots, \ell$.
We collect these bounds by setting $g\in\cD=[g_1^-, \, g_1^+]\times \dots \times [g_{\ell}^-, \, g_{\ell}^+]\subset\R_+^{\ell}$, where the parameter set $\cD$ contains all restricting intervals.
However, it is also possible that no a priori knowledge of the parameter set is given, then we define $\cD=\R^{\ell}_+$.

Our goal is to design the damping values based on some optimization of an objective function.
The problem of finding optimal external dampers was widely investigated in the literature, see \cite{Kan13,MueS85,Tak97,Gen09,Inm06,Ves11}.
The objective of this work, is to solve the problem of semi-active damping, i.e., for a given vibrational system of the form \eqref{eq:SOsys}, to determine the best damping $D(g)$ that ensures optimal attenuation of the output $y$.
Since we want to minimize the maximum response gain for the entire range of time, we minimize the $L_\infty$-norm of the output $y$.
For that, we can make use of the following bound 
\[
\|y\|_{L_{\infty}}\leq  \|\cG(\cdot, g)\|_{\cH_2}\| u\|_{L_2},
\]
where $\cG(s,g):=C(s^2M + sD(g) + K)^{-1}B$ is the \emph{transfer function} associated to the viabrational system in \eqref{eq:SOsys} that describes input to output behavior in the frequency domain. 
We see that the norm of $y$ is bounded by the $\cH_2$-norm of the transfer function $\cG(\cdot, g)$ that is also called the \emph{energy response} and defined as 
\begin{align*}
\bJ(g):=\|\mathcal{G}(\cdot;g)\|_{\mathcal{H}_2} = \left( \frac{1}{2\pi} \int_{-\infty}^{\infty} \mathrm{trace}(\mathcal{G}(\mathrm{i}w;g)^{\rmH}\mathcal{G}(\mathrm{i}w;g))\mathrm{d}w \right)^{\frac{1}{2}}.
\end{align*}
In this article, we will optimize the damping values in such a way that the $\cH_2$-norm of the transfer function, and hence the energy response, is minimized.
This criterion was also used in \cite{morTomBG18,morGugAB08,morBenKTetal16}.

The vibrational system \eqref{eq:SOsys}  can be represented as a first order system
\begin{align}\label{eq:FOsys}
\begin{split}
\dot{z}(t) &= \cA(g)z(t) + \cB u(t),\\
y(t) &= \cC z(t)
\end{split}
\end{align}
with
\[
\cA(g):=\begin{bmatrix}
0 & I \\
-M^{-1}K & -M^{-1}D(g)
\end{bmatrix}, \quad \cB := \begin{bmatrix}
0\\
M^{-1}B
\end{bmatrix}, \quad \cC := \begin{bmatrix}
C & 0
\end{bmatrix}\;\text{  and  }\; z(t) := \begin{bmatrix}
x(t)\\
\dot{x}(t)
\end{bmatrix}.
\]
The latter is a \emph{linearization} of \eqref{eq:SOsys}, and we emphasize that there are many other useful linearizations that could be used here. We use \eqref{eq:FOsys} for simplicity.
As described in \cite{ZhoDG96}, the $\cH_2$-norm of the transfer function can be computed as 
\begin{align}\label{eq:SysResTr}
\bJ(g)=\trace{\cC \cP(g) \cC^\T}^{\frac{1}{2}} 
= \trace{C P_{11}(g) C^\T}^{\frac{1}{2}},
\end{align}
where the matrix $\cP(g)$ is called \emph{controllability Gramian} and is defined in the time domain or the frequency domain as 
\begin{align}\label{eq:Gramian}
\begin{split}
\cP(g) = \begin{bmatrix}
P_{11}(g) & P_{12}(g) \\ P_{12}(g)^{\T} & P_{22}(g)
\end{bmatrix}:=&{} \int_{0}^{\infty}e^{\cA(g) t}\cB\cB^{\T}e^{\cA(g)^{\T} t}\mathrm{d} t\\
:=&{} \frac{1}{2\pi}\int_{-\infty}^{\infty} (i\omega I - \cA(g))^{-1}\cB\cB^{\T}(-i\omega I - \cA(g))^{-\T}\mathrm{d}\omega.
\end{split}
\end{align}
The upper left block $P_{11}(g)$ of the Gramian $\cP(g)$ is called \emph{position controllability Gramian} and encodes the reachabillity space of the state $x(t)$ of the corresponding second-order system \eqref{eq:SOsys}.
The controllability Gramian $\cP(g)$ is computed by solving the Lyapunov equation 
\begin{align}\label{eq:Lyap}
\cA(g)\cP(g) + \cP(g)\cA(g)^\T = -\cB\cB^\T.
\end{align}

In order to find the damping gains $g\in\cD$ that minimize the energy response $\bJ(g)$, we have to solve a Lyapunov equation \eqref{eq:Lyap} in every step of the optimization method.
Since the Lyapunov equation solves are very demanding if the matrices are of large dimensions, the minimization process would lead to high computational cost and hence be inefficient or unfeasible in a large-scale setup.
In the literature, several approaches have been developed to accelerate the minimization process.
The authors in \cite{morBenKTetal16} utilize the dominant pole algorithm in order to build a reduced minimization problem that is fast solvable.
On the other hand, in \cite{morTomBG18}, an efficient optimization approach using structure-preserving parametric
model reduction based on the iterative rational Krylov algorithm (\texttt{sym2IRKA}) is used to derive an efficient optimization algorithm.
In \cite{morBeaGT20}, a sampling-free approach is presented that reduces the system \eqref{eq:SOsys} for all admissible parameters.
Alternatively, in \cite{morTomV20} the $\cH_{\infty}$-norm of the transfer function $\cG$ is minimized.
In \cite{morBenTT11a,morTru04,morTruV09}, the authors present different reduction techniques to optimize the related problem of minimizing the total average energy for the system \eqref{eq:SOsys} with no input.
It is worth mentioning that the problem of optimizing damping positions is investigated in \cite{KanPTetal19,morTom11,morBenTT11,GurM92,DymD21}. 
This  is a challenging problem especially for large-scale systems and it is not the focus of this work.

In this paper, we apply the reduced basis method (RBM) and modifications of it, in order to reduce the Lyapunov equation \eqref{eq:Lyap} such that the surrogate equation can be solved in a reasonable time.
The RBM is a well-established method to reduce parameter-dependent partial differential equations \cite{morHesRS16,morQuaMN16,morVerPRetal03,morVerPP03,morVerP05}.
Later, it was used for Riccati equations \cite{morSchH15} and finally the RBM was applied to Lyapunov equations by Son an Stykel in \cite{morSonS17}. 
In \cite{morPrzV21}, the authors use the RBM in order to reduce parametric differential-algebraic systems. 
The method is decomposed into an offline and an online phase.
In the \emph{offline phase} of the RBM, a reduced space is determined that approximates the solution space of the original Lyapunov equation for all parameters $g\in\cD$.
To evaluate the quality of the approximated reduced space, we suggest two error estimators that evaluate the error in the controllability Gramians and in the corresponding energy responses for different damping values.
The approximated solution space is then used in the \emph{online phase} to solve the corresponding reduced Lyaponv equations, leading to an approximation of the energy response determined by Equation \eqref{eq:SysResTr} for all requested parameters $g\in\cD$.
To avoid adding informations into the basis that are not needed during the optimization, we introduce an adaptive method that enriches the basis within the optimization process.

The rest of this paper is structured as follows: In \Cref{sec:RBM}, we describe the reduced basis method, and introduce error estimators that are suitable for our method.
The adaptive enrichment method is then introduced in \Cref{sec:Hybrid}.
Afterwards, we describe some implementation details in \Cref{sec:ImplDetail}.
In order to illustrate the effectiveness of our methods, in \Cref{sec:Results}, we apply our methods to two examples, and finally \Cref{sec:Conclusion} concludes the manuscript.


\section{Reduced basis method}\label{sec:RBM}
As described above, we want to find the damping values that minimize the energy response $\bJ(g)$. The optimization procedure requires the evaluation of this energy  for several parameters $g\in\cD$. 
To accelerate the optimization, we want to speed up the computation of the Gramian $\cP(g)$ or $P_{11}(g)$, which contains the numerically most costly step. The reduced basis method (RBM) has been found to be an effective tool to deal with this problem, see \cite{morSonS17}. 
We apply the RBM in order to reduce the optimization problem in \Cref{ssec:Reduction}.
Afterwards, in \Cref{ssec:ErrEst}, we derive error estimators that are needed in the RBM to evaluate the quality of the approximations.

\subsection{Reduction using the reduced basis method}\label{ssec:Reduction}

In order to accelerate the computation of the position controllability Gramians $P_{11}(g)$, we aim to find a space $\cV_{\cD}$ that spans the controllability space of the second-order system \eqref{eq:SOsys} and hence the columns of the Gramians $P_{11}(g)$, i.e.
\[
\myspan{P_{11}(g)}\subset\cV_{\cD}=\myspan{V_{\cD}},
\]
for all admissible parameters $g\in\cD$, where $V_{\cD}$ is a basis that spans the space $\cV_{\cD}$.
Then, for all $g\in\cD$, there exists a matrix $X_{11}(g)$ with $P_{11}(g)=V_{\cD}X_{11}(g)V_{\cD}^{\T}$.
Since such a space is in general not available, the idea of the RBM is to find a reduced space $\cV=\myspan{V_1}$, $V_1\in\Rnr$ so that the position controllability Gramian $P_{11}(g)$ can be approximated as 
\begin{align}\label{eq:P_approx}
P_{11}(g) \approx \widetilde{P}_{11}(g) = V_1\widehat{P}_{11}(g)V_1^{\T}
\end{align}
for suitable $\widehat{P}_{11}(g)\in\Rrr$, and for all damping parameters $g\in\cD$.
The matrices $\widehat{P}_{11}(g)$ are determined by solving Lyapunov equations of dimension $2r$, which are fast computable if $r$ is small enough.

The RBM is decomposed into an offline and an online phase.
The offline phase consists of the computation of the basis $V_1\in\Rnr$, which is rather time consuming but needs to be performed only once.
In the online phase, on the other hand, Lyapunov equations are solved on the reduced space
$\cV$, which is fast and can be performed multiple times.
This phase involves performing the optimization, which requires evaluating the approximated energy response for multiple damping parameters.
In order to describe the two phases in more detail, we need a criterion to evaluate the quality of the reduced space $\cV$.
Thus, we assume that we have an error estimator $\Delta(g)$ that provides a criterion to determine how well the solution space for a parameter $g$ is approximated by the current basis $V_1$.
The used error estimators are described later in \Cref{ssec:ErrEst}.

\emph{Offline phase}: 
In this phase, we aim to find an approximation $\cV$ of the space $\cV_{\cD}$ and the corresponding basis $V_1$, which will be used later to define an approximation of the position controllability Gramians $P_{11}(g)$ for the required parameters $g\in\cD$.
Since we can not evaluate and investigate an infinite number of parameters that are given in $\cD$, we define a test-parameter set $\cD_{\Test}\subset\cD$ that is finite and well distributed in $\cD$. 
Finding a good sampling strategy for the test-parameter set is a challenging task, which was addressed in \cite{morCheFB21,morRozHP08,morHesSZ14,morEftKP11} for partial differential equations.
In the following, we derive a space $\cV$ that includes approximations of the solution spaces corresponding to all test-parameters $g\in\cD_{\Test}$.
Then, if $\cD_\Test$ is well-chosen, the space $\cV$ approximates the solution space for all parameters in $\cD$.

Because of the low-rank structure of the right-hand side in \eqref{eq:Lyap}, we can assume that the solution $\cP(g)$ of the Lyapunov equation in \eqref{eq:Lyap} is well approximated by a low-rank factor $\cZ(g)$ such that 
\begin{align}\label{eq:LRStruc}
\cP(g)\approx \cZ(g)\cZ(g)^{\T},\qquad P_{11}(g)\approx Z_1(g)Z_1(g)^{\T},\qquad\text{where}\qquad \cZ(g)= \begin{bmatrix}
Z_1(g) \\ Z_2(g)
\end{bmatrix},
\end{align}
for $g\in\cD$.
In order to build a first basis $V_1$, we choose one arbitrary test-parameter $g_0\in\cD_\Test$ and solve the Lyapunov equation in \eqref{eq:Lyap} to obtain the low-rank factor $\cZ(g_0)$ that includes $Z_1(g_0)$. 
Using the low-rank factor $Z_1(g_0)$, we define the basis
\[
V_1:=\orth(Z_1(g_0)).
\]
The operator $\orth(\cdot)$ describes the orthonormalization of the columns of a given matrix.  
After forming our first basis, we evaluate the quality of the approximation of the position controllability Gramian $P_{11}(g)$ for all remaining parameters $g\in\cD_{\Test}$.
To this aim, we compute the error estimate $\Delta(g)$ for all these parameters and define the largest error estimate as
\[
\Delta^{\max} := \Delta(g_{1}) := \max_{g\in\cD_\Test}\Delta(g),
\]
where $g_{1}$ is the parameter that leads to the largest estimate.
If $\Delta^{\max}$ is larger than a given tolerance $\tol_f$, we know that the current basis does not approximate the solution space good enough for at least one parameter $g_{1}$.
Hence, we need to enlarge the basis $V_1$. 
An obvious candidate to enrich the basis are the columns of the Gramian $P_{11}(g_1)$ for the parameter $g_{1}$ that results in this largest error estimate.
We compute the low-rank factor $Z_1(g_{1})$ by solving the Lyapunov equation in \eqref{eq:Lyap} and set \[V_1 = \orth([V_1, Z_1(g_{1})])\] to build the new basis.
We continue with this procedure until the error estimate $\Delta(g)$ is smaller than the tolerance $\tol_f$ for every damping value $g \in \cD_{\Test}$.

This method is described in \Cref{algo:RBMOff}.
In Step \ref{StepM1} and \ref{StepM2} we define the set of already used parameters $\cM$ so that we do not evaluate their error estimates in Step \ref{Sr5} and \ref{S11}.
\begin{algorithm}
	\caption{Offline Phase - RBM}\label{AlgoOfflineRBM}
	\label{algo:RBMOff}
	\textbf{Input:} $\cA:\cD\to\Rnn$ asymptotically stable, $\cB\in\Rnm$, test-parameter set $\cD_{\Test}$, tolerance $\tol_f$\\
	\textbf{Output:} Orthonormal basis $V_1$
	\begin{spacing}{1.14}
		\begin{algorithmic}[1] 
			\State Choose any $g_{0}\in\cD_{\Test}$.
			\State Solve the Lyapunov equation in \eqref{eq:Lyap} at $g_{0}$ to obtain $Z_1(g_{0})$.
			\State Set $\cM:=\{g_{0}\}$.\label{StepM1}
			\State Set $V_1 := \mathrm{orth}(Z_1(g_{0}))$.\label{StepONB1}
			\State Set $k=1$.
			\State Determine $g_1 := \argmax_{g\in\cD_{\Test}\setminus\cM}\Delta(g)$.\label{Sr5}
			\State Set $\Delta^{\max}=\Delta(g_1)$.
			\While{$\Delta^{\max}>\tol_f$}
			\State Solve Lyapunov equation in \eqref{eq:Lyap} at $g_k$ to obtain $Z_1(g_k)$.
			\State Set $\cM := \cM\cup\{g_k\}$.\label{StepM2}
			\State Set $V_1 := \mathrm{orth}([V_1,~Z_1(g_k)])$.\label{StepONB2}
			\State Determine $g_{k+1}:= \argmax_{g\in\cD_{\Test}\setminus\cM}\Delta(g)$.\label{S11}
			\State Set $\Delta^{\max}:=\Delta(g_{k+1})$.
			\State Set $k = k+1$.
			\EndWhile
		\end{algorithmic}
	\end{spacing}
\end{algorithm}

\newpage

\emph{Online phase}: 
After the basis $V_1$ is computed, it can be used to determine an approximation of the position controllability Gramian $P_{11}(g)$ as described in \eqref{eq:P_approx}. 
For that, define the basis 
\[
V := \begin{bmatrix}
V_1 & 0 \\ 0 & V_1
\end{bmatrix}
\] 
corresponding to the first-order system in \eqref{eq:FOsys} and define the by $V$ reduced matrices 
\[
\widehat{\cA}(g):= V^{\T}\cA(g)V
,\qquad \widehat{\cB}:=V^{\T}\cB\quad\text{ and }\quad \widehat{\cC}:=\cC V.
\]
We compute the reduced Gramian $\widehat{P}_{11}(g)$ that is used in \eqref{eq:P_approx} by solving the reduced Lyapunov equation 
\begin{align}\label{eq:redLE}
\widehat{\cA} \widehat{\cP}(g) +  \widehat{\cP}(g)\widehat{\cA}^{\T} = -\widehat{\cB}\widehat{\cB}^{\T},
\end{align}
where $\widehat{\cP}(g) = \begin{bsmallmatrix}
\widehat{P}_{11}(g) & \widehat{P}_{12}(g) \\ \widehat{P}_{12}(g)^{\T} & \widehat{P}_{22}(g)
\end{bsmallmatrix} = \widehat{\cZ}(g) \widehat{\cZ}(g)^{\T}$ and $\widehat{\cZ}(g) = \begin{bsmallmatrix}
\widehat{Z}_1(g) \\ \widehat{Z}_2(g)
\end{bsmallmatrix}$. 
Note that the projected Lyapunov equation in \eqref{eq:redLE} is of dimension $2r$ and can therefore be solved cheaply compared to the original Lyapunov equation in \eqref{eq:Lyap}.

In practice, we make use of the low-rank structure of the Gramians described in \eqref{eq:LRStruc} so that we obtain the approximations
\[
Z_1(g)\approx \tZ_1(g) = V_1 \widehat{Z}_1(g)\qquad\text{or}\qquad 
\cZ(g)\approx \widetilde{\cZ}(g) = V \widehat{\cZ}(g).
\]

\emph{Optimization}: 
Finally, we include the online phase of the RBM into our optimization problem.
We assume that the bases $V_1$ and $V$ were computed in the offline phase, beforehand.
Then, in any step of the optimization process (in a parameter $g$), we solve the reduced Lyapunov equation from \eqref{eq:redLE} so that the reduced Gramians $\widehat{\cP}(g)$ and $\widehat{P}_{11}(g)$ are available.
Afterwards, we define the approximated energy response 
\begin{align}\label{eq:Min_Red}
\tilde{\bJ}(g) := \mathrm{tr}(C V_1\widehat{P}_{11}(g)V_1^{\T}C^{\T})^{\frac{1}{2}} = \mathrm{tr}(\cC V \widehat{\cP}(g)V^{\T}\cC^{\T})^{\frac{1}{2}}
\end{align}
that is minimized.
The computation of the reduced energy response in \eqref{eq:Min_Red} includes the solving of a Lyapunov equation from \eqref{eq:redLE} of dimension $2r$ which is performed in every step of the optimization process.
This is a significant acceleration compared to the optimization of the original energy response in \eqref{eq:SysResTr} where the Lyapunov equation from \eqref{eq:Lyap} of dimension $2n$ needs to be solved in every iteration step.

\subsection{Error estimator}\label{ssec:ErrEst}
For the reduced basis method presented above, as well as for the adaptive method presented in the following section, error estimators are needed to evaluate the quality of the resulting approximations.
In this section, we will derive error estimators of two different quantities.

The first one is the approximation error of the position controllability Gramians, defined as 
\[
\|\cE_{11}(g)\| := \|P_{11}(g) - \widetilde{P}_{11}(g)\|,
\] 
where $\widetilde{P}_{11}(g)$ is an approximation of $P_{11}(g)$.
The \emph{position controllability Gramian error} $\cE_{11}(g)$ is the upper left block of the corresponding error
\[
\cE(g)=\begin{bmatrix}
\cE_{11}(g) & \cE_{12}(g)\\ \cE_{12}(g)^{\T} & \cE_{22}(g)
\end{bmatrix} = \begin{bmatrix}
P_{11}(g) & P_{12}(g)\\ P_{12}(g)^{\T} & P_{22}(g)
\end{bmatrix} -\begin{bmatrix}
\widetilde{P}_{11}(g) & \widetilde{P}_{12}(g)\\ \widetilde{P}(g)^{\T} & \widetilde{P}_{22}(g)
\end{bmatrix} = \cP(g)-\widetilde{\cP}(g).
\]
In the literature there are various upper bounds for $\|\cE(g)\|$, i.e. for the error in the solution of the Lyapunov equation in \eqref{eq:Lyap}, see e.g. \cite{morSonS17,morPrzV21}.
These bounds use the corresponding residual
\begin{equation}\label{eq:ResLE}
\cR(g):=\cB\cB^\T+\cA(g)\widetilde{\cP}(g)+\widetilde{\cP}(g)\cA(g)^\T
\end{equation}
divided by the coercivity constant of the Lyapunov equation in \eqref{eq:Lyap} which is the minimal eigenvalue of the corresponding linear system matrix. 
For vibrational systems with small internal damping, the smallest eigenvalues are close to the imaginary axis, and hence, the coercivity constant is small so that the error estimators lead to large, very conservative values and are not feasible for our application.
Additionally, we only want to evaluate the approximation of the controllability space of $x(t)$ encoded in $P_{11}(g)$, while evaluating the approximation of the complete Gramian $\cP(g)$ contains information that are not used.

The second quantity that we consider is the approximation error of the function values of the energy response, i.e.
\begin{align*}
\cE_{\bJ}(g):=\left|\bJ(g) - \widetilde{\bJ}(g)\right| = \left|\trace{\cC \cP(g)\cC^{\T}} - \trace{\cC\widetilde{\cP}(g)\cC^{\T}}\right|
= \left|\trace{\cC\cE(g)\cC^{\T}}\right|= \left|\trace{C\cE_{11}(g)C^{\T}}\right|.
\end{align*}

We notice that the error $\cE(g)$ or $\cE_{11}(g)$ is required for both, the evaluation of the error $\|\cE_{11}(g)\|$ as well as for the computation of the error $\cE_{\bJ}(g)$.
Since it is costly to compute the error $\cE(g)$ for every requested parameter $g$, we require an approximate $\widetilde{\cE}(g)\approx\cE(g)$.
We aim to find such an $\widetilde{\cE}(g)$ in the following.
Firstly, notice that the error $\cE(g)$ is the solution of the following Lyapunov equation
\begin{align}\label{eq:ErrEq}
\cA(g)\cE(g)+ \cE(g)\cA(g)^\T = -\cR(g)
\end{align}
that is also called \emph{error equation}.
Hence, a second RBM that follows the same scheme as presented in \Cref{ssec:Reduction} can be applied to determine a basis $V_{1,\err}$ and to approximate the error as 
\begin{align*}
\cE_{11}(g) \approx \widetilde{\cE}_{11}(g) = V_{1,\err}\widehat{\cE}_{11}(g)V_{1,\err}^{\T}.
\end{align*}
To avoid confusion, we denote the second reduced basis method that determines the basis $V_{1,\err}$ for the error estimation as EE-RBM.
The basis $V_{1,\err}$ can be constructed using the solutions of \eqref{eq:ErrEq} for some parameters.
However, the right-hand side $\cR(g)$ of the error equation in \eqref{eq:ErrEq} does not consist of low-rank factors, and hence, solving this Lyapunov equation is generally numerically costly, since solution algorithms cannot exploit such a low-rank factor structure.
Hence, in order to avoid having do solve the error equation in \eqref{eq:ErrEq}, we investigate the error $\cE_{11}(g)$ and the structure of the corresponding error spaces in more detail to simplify its computation.
For this purpose, we write the position controllability Gramian as $P_{11}(g)=V_{\cD}X_{11}(g)V_{\cD}^{\T}$, where $V_{\cD}$ is a basis spanning the (complete) solution space $\cV_{\cD}$ of the second-order system in \eqref{eq:SOsys} for all parameters $g\in\cD$.
The error is then given as 
\[
\cE_{11}(g) = V_{\cD}X_{11}(g)V_{\cD}^{\T} - V_1\widehat{P}_{11}(g)V_1^{\T}\] 
and hence we obtain, that the error $\cE_{11}(g)$ lies in the space spanned by the basis $V_{\cE} = \orth([V_1, V_{\cD}])$ for all parameters, which was investigated in \cite{morCheFdetal21} for linear systems.
Since $V_1$ is known already from the first RBM, the remaining task is to determine $V_{\cD}$.
However, obviously, the basis $V_{\cD}$ is not available, otherwise we would have a basis that spans the solution space of the Lyapunov equation in \eqref{eq:Lyap} for all parameters without any error.
This is why we apply the second reduced basis method (EE-RBM) and derive an approximation of $V_{\cE}$ that is called $V_{1,\err}$.
However, because of the structure of the basis $V_{\cE}$, we can solve a second Lyapunov equation from  \eqref{eq:Lyap} instead of the error equation in \eqref{eq:ErrEq}, which is of a more advantageous structure because of the low-rank factor structure of the right-hand side.
Hence, in every step of EE-RBM, we solve the Lyapunov equation from \eqref{eq:Lyap} in a certain parameter $g^{\rr}$ to obtain the corresponding solution $Z_1(g^{\rr})$ and to enrich the basis of the error equation in \eqref{eq:ErrEq} as $V_{1,\err} = \orth([V_{1,\err},~ V_1,~ Z_1(g^{\rr})])$, and therefore build a basis that approximates $V_{\cE}$.
Adding $V_1$ and $Z_1(g^{\rr})$ is equivalent to solving the error equation from \eqref{eq:ErrEq} in $g^{\rr}$ and adding the resulting solution to the basis $V_{1,\err}$. 

In order to include the computation of the basis $V_{1,\err}$ into the first RBM, we run both methods, the RBM and the EE-RBM in parallel.
The first parameters $g_0$ and $g_0^{\rr}$ are chosen arbitrarily in $\cD_{\Test}$ with $g_0^{\rr}\neq g_0$.
We refer to \cite{morFenB21} for a mathematical analysis why $g_0^{\rr}$  should be chosen different from $g_0$.
We compute the basis $V_1=\orth\left(Z_1(g_0)\right)$ as described in the previous subsection and, in addition, solve the Lyapunov equation from \eqref{eq:Lyap} in $g_0^{\rr}$ to obtain $Z_1(g_0^{\rr})$ such that our first error space basis is given as 
\[
V_{1,\err}=\orth\left([V_1, Z_1(g_0^{\rr})]\right)=\orth\left([Z_1(g_0), Z_1(g_0^{\rr})]\right).
\]
As described above, the next parameter $g_1$ is the one that leads to the largest error estimate $\Delta(g)$ and we use the corresponding solution $Z_1(g_1)$ to enrich the basis $V_1$.
Next, the parameter $g_1^{\rr}$ is chosen to be that one that results in the largest residual of the error equation in the Frobenius norm, i.e. the parameter $g_1^{\rr}$ that leads to the largest value
\begin{equation}\label{eq:Resisual}
\|R_{\rr}(g)\|_{\Fr} := \|\cA(g)\widetilde{\cE}(g)+ \widetilde{\cE}(g)\cA(g)^\T + \cB\cB^\T +\cA(g)\widetilde{P}(g) + \widetilde{P}(g)\cA(g)^\T\|_{\Fr} .
\end{equation}
We use the parameters $g_1^{\rr}$ to generate the next error equation basis 
\[
V_{1,\err}=\orth\left([V_{1,\err}, ~V_1, ~ Z_1(g_1^{\rr})]\right) =\orth\left([V_{1,\err},~ Z_1(g_1),~ Z_1(g_1^{\rr})]\right).
\]
We continue with this process until the maximum error estimator $\Delta(g)$ is smaller than a certain tolerance for all parameters and the overall RBM is finished.

After we have determined an error equation basis $V_{1,\err}$, this basis is used to derive the corresponding error estimators. 
For that, we define 
\[
V_{\err}:= \begin{bmatrix}
V_{1,\err} & 0 \\ 0 & V_{1,\err}
\end{bmatrix},
\]
where $V_{\err}$ spans an approximation of the solution space of the error equation in \eqref{eq:ErrEq}.
Using this basis, we define the reduced matrices 
$
\widehat{\cA}_{\err}(g) := V_{\err}^{\T}\cA(g)V_{\err}$ and $\widehat{\cR}_{\err}(g) :=
V_{\err}^{\T}\cR(g)V_{\err}
$
that lead to  the reduced error equation
\begin{align}\label{eq:redErrEq}
\widehat{\cA}_{\err}(g)\widehat{\cE}(g)+ \widehat{\cE}(g)\widehat{\cA}_{\err}(g)^\T  = -\widehat{\cR}_{\err}(g)
\end{align}
that is solved by the reduced error matrix $\widehat{\cE}(g)$, similar to the step in \eqref{eq:redLE} for the RBM.
The approximation of the position controllability Gramian error $\widetilde{\cE}_{11}(g)$ is then built as 
\begin{align}\label{eq:redErr}
\widetilde{\cE}(g) = \begin{bmatrix}
\widetilde{\cE}_{11}(g) & \widetilde{\cE}_{12}(g)\\
\widetilde{\cE}_{12}(g)^{\T} & \widetilde{\cE}_{22}(g)
\end{bmatrix} 
=\begin{bmatrix}
V_{1,\err}\widehat{\cE}_{11}(g)V_{1,\err}^{\T} & V_{1,\err}\widehat{\cE}_{12}(g)V_{1,\err}^{\T}\\
V_{1,\err}\widehat{\cE}_{12}(g)^{\T}V_{1,\err}^{\T} & V_{1,\err}\widehat{\cE}_{22}(g)V_{1,\err}^{\T}
\end{bmatrix} 
= V_{\err}\widehat{\cE}(g)V_{\err}^{\T}.
\end{align}

Using the error approximation $\widetilde{\cE}(g)$ and $\widetilde{\cE}_{11}(g)$ from \eqref{eq:redErr}, we define the two error estimators
\begin{align}\label{eq:ErrEst}
\Delta_{1}(g)
:=\left|\trace{\cC\widetilde{\cE}(g)\cC^{\T}}\right|
=\left|\trace{C\widetilde{\cE}_{11}(g)C^{\T}}\right|
\qquad\text{and}\qquad
\Delta_{2}(g)
:=\left\|\widetilde{\cE}_{11}(g)\right\|_{\mathrm{F}}.
\end{align}

We use the expression $\Delta(g)$ in the following and leave it to the user to choose the more appropriate one.
The first RBM combined with the EE-RBM results in \Cref{algo:ErrEst}.
\begin{algorithm}
	\caption{Error estimation within the RBM}\label{algo:ErrEst}
	\textbf{Input:} $\cA:\cD\to\Rnn$ asymptotically stable, $\cB\in\Rnm$, test-parameter set $\cD_{\Test}$, tolerance $\tol_f$\\
	\textbf{Output:} Orthonormal bases $V_1$, $V_{1,\err}$
	\begin{spacing}{1.14}
		\begin{algorithmic}[1] 
			\State Choose any $g_0,~g_0^{\rr}\in\cD_{\Test}, ~ g_0\neq g_0^{\rr}$.
			\State Solve the Lyapunov equation \eqref{eq:Lyap} at $g_{0}$ to obtain $Z_1(g_{0})$.
			\State Set $\cM:=\{g_{0}\}$.
			\State Set $V_1 := \mathrm{orth}(Z_1(g_{0}))$.
			\State Solve Lyapunov equation \eqref{eq:Lyap} for $g_0^{\rr}$ to obtain $Z_1(g_0^{\rr})$.
			\State Set $V_{1,\err} := \orth([Z_1(g_{0}), ~Z_1(g_0^{\rr})])$.
			\State Set $k:=1$.
			\State Determine $g_1 := \argmax_{g\in\cD_{\Test}\setminus\cM}\Delta(g)$.
			\State Set $\Delta^{\max} := \Delta(g_1)$.
			\State Determine $g_{1}^{\rr}:= \argmax_{g\in\cD_{\Test}\setminus\cM}\|\cR_{\rr}(g)\|_{\Fr}$.
			\While{$\Delta^{\max}>\tol_f$}
			\State Solve Lyapunov equation \eqref{eq:Lyap} at $g_k$ to obtain $Z_1(g_k)$.
			\State Set $\cM := \cM\cup\{g_k\}$.
			\State Set $V_1 := \mathrm{orth}([V_1,~Z_1(g_k)])$.
			\State Solve Lyapunov equation \eqref{eq:Lyap} for $g_k^{\rr}$ to obtain $Z_1(g_k^{\rr})$.
			\State Set $V_{\err} := \orth([V_{1,\err},~ Z_1(g_k), ~Z_1(g_k^{\rr})])$.
			\State Determine $g_{k+1}:= \argmax_{g\in\cD_{\Test}\setminus\cM}\Delta(g)$.
			\State Set $\Delta^{\max}:=\Delta(g_{k+1})$.
			\State Determine $g_{k+1}^{\rr}:= \argmax_{g\in\cD_{\Test}\setminus\cM}\|\cR_{\rr}(g)\|_{\Fr}$.
			\State Set $k := k+1$.
			\EndWhile

		\end{algorithmic}
	\end{spacing}
\end{algorithm}
We observe that the first steps of the RBM with the EE-RBM lead to rough error estimates since the basis $V_{1,\err}$ includes only a few solutions. 
However, the larger and therefore better the basis $V_1$ is, the larger and more detailed is the basis $V_{1,\err}$.

\begin{remark}
Numerical experiments suggest that adding the solution vectors of $Z_1(0)$ corresponding to the undamped system to our basis $V_1$ leads to more robust results.
Since this basis is independent of the damping values, we compute it beforehand and initialize the basis $V_1 = \orth(Z_1(0))$.
\end{remark}

\begin{remark}
In practice, the computation of $\cR_{\rr}(g)$ can be performed efficiently as follows. We make use of the trace formulation of the Frobenius norm and
utilize the low-rank representations $\widetilde{\cE}(g) = V_{\err}\widehat{\cE}(g)V_{\err}^{\T}$ and $\widetilde{P}(g)=V \widehat{P}(g)V^{\T}$ and the trace properties to obtain the fast computable residual representation:
\begin{align*}
\|\cR_{\rr}(g)\|_{\Fr}^2 &= 2\trace{V_{\err}^{\T}\cA(g)V_{\err}\widehat{\cE}(g)V_{\err}^{\T}\cA(g)V_{\err}\widehat{\cE}(g)} 
+ 2\trace{V_{\err}^{\T}\cA(g)^{\T}\cA(g)V_{\err}\widehat{\cE}(g)V_{\err}^{\T}V_{\err}\widehat{\cE}(g)} \\
&+ 4\trace{V_{\err}^{\T}(g)\cB\cB^{\T}\cA(g)V_{\err}\widehat{\cE}(g)}
+ 2\trace{V^{\T}\cA(g)V \widehat{P}(g)V^{\T}\cA(g)V \widehat{P}(g)} \\
&+ 2\trace{V^{\T}\cA(g)^{\T}\cA(g)V \widehat{P}(g)V^{\T}V \widehat{P}(g)} 
+ 4\trace{V^{\T}\cB\cB^{\T}\cA(g)V \widehat{P}(g)}\\
&+ 4\trace{V^{\T}\cA(g)V_{\err}\widehat{\cE}(g)V_{\err}^{\T}\cA(g)V \widehat{P}(g)} 
+ 2\trace{V^{\T}V_{\err}\widehat{\cE}(g)V_{\err}^{\T}\cA(g)^{\T}\cA(g)V \widehat{P}(g)} \\
&+ \trace{\cB^{\T}\cB\cB^{\T}\cB}.
\end{align*}

\end{remark}


\section{Adaptive basis building}\label{sec:Hybrid}
In the previous section, we proposed an RBM leading to a space $\cV$ that approximates the controllability space of the state $x(t)$ of system \eqref{eq:SOsys} for all parameters in $\cD$. 
This space is then used to derive a reduced optimization problem. 
However, if the parameter set is equal to $\cD = \R_+^{\ell}$, i.e., there is no prior knowledge about the parameter set, the method is not applicable.
Additionally, the optimization process might only use parameters from a subset of $\cD$ such that the basis $V_1$ of the space $\cV$ from \Cref{sec:RBM} contains controllability space information corresponding to parameters $g$ that are not used throughout the optimization process and therefore the basis $V_1$ might be of too large dimension. 
This motivates an adaptive scheme, i.e. a procedure, that enriches the basis $V_1$ within the optimization process based on the quality of the approximation in the considered parameters. 
To find out whether the current basis $V_1$ is adequate for the current parameter $g$, we again need an error estimator $\Delta(g)$. 
In what follows, we will firstly describe this adaptive RBM in \Cref{ssec:ReductionARBM}  assuming an error estimator is provided. Then, in \Cref{ssec:ErrEstARBM}, an error estimator is proposed for this adaptive procedure.

\subsection{Reduction using the adaptive reduced basis method}\label{ssec:ReductionARBM}

The idea of the adaptive RBM is to enrich the basis $V_1$ within the optimization process. Consequently, there is no offline phase for this methodology, since the reduced bases are constructed within the optimization steps.
Firstly, we select a parameter $g_0\in\cD$ as the initial value for the optimization process and solve the Lyapunov equation in \eqref{eq:Lyap} for this parameter to obtain the low-rank factor $Z_1(g_0)$.
We set the first basis to be $V_1=\orth(Z_1(g_0))$ that is used to define the reduced optimization problem in \eqref{eq:Min_Red}, which depends on the solution of the reduced Lyapunov equation in \eqref{eq:redLE}.
In contrast to the previous method, we add an additional stopping criterion within the optimization process that interrupts the procedure whenever the solution space corresponding to the current parameter $g$ is not well-approximated by the basis $V_1$.
To achieve this stopping, we modify the goal function as described by \Cref{algo:hybridRedOpt}.
In every iteration of the minimization, we query the error estimate $\Delta(g)$ of the current parameter $g$ as described in Step \ref{S2}.
If the error estimate is smaller than a given tolerance, we proceed with the function evaluation in Steps \ref{Sf1} and \ref{Sf2} to obtain the resulting function value $\widetilde{\bJ}(g)$ and continue with the minimization.
On the other hand, if the error estimate is larger than the tolerance, this means that the current basis $V_1$ does not approximate the solution space of the Lyapunov equation \eqref{eq:Lyap} for the current parameter $g$ sufficiently well.
 Hence, we return that the minimization did not converge.
In this case, we need to enrich the basis $V_1$.
Therefore, we solve the Lyapunov equation from \eqref{eq:Lyap} in this parameter $g$ to obtain $Z_1(g)$ and define the updated basis 
\[V_1=\orth([V_1, Z(g)]).\]
Consequently, we obtain a new optimization problem \eqref{eq:Min_Red} that is defined by the new basis $V_1$ and the new projected Lyapunov equation from \eqref{eq:redLE}.
Since the optimized function depends on the current basis $V_1$, which changes during the optimization procedure, convergence problems may occur.
In this case, we restart the optimization procedure after we enrich the basis and use the current parameter $g$ as initial value.
We continue with this procedure until the optimum is reached.

\subsection{Error Estimator}\label{ssec:ErrEstARBM}

Finally, we have a closer look at the error estimator $\Delta(g)$ for the adaptive RBM.
We follow the same idea as in \Cref{ssec:ErrEst} and run a second reduced basis method to generate a basis $V_{1,\err}$ that spans an approximation of the error space.
The equation in \eqref{eq:ErrEst} defines then the error estimators $\Delta_1(g)$ and $\Delta_2(g)$ corresponding to the bases $V_1$ and $V_{1,\err}$.
In this adaptive procedure, the basis $V_{1,\err}$ is enlarged whenever the basis $V_1$ is expanded.
In this way, the error approximation, and thus the error estimator, becomes more accurate the closer we get to the optimizing parameter.

The detailed procedure is described in \Cref{algo:adaptiveRBM}.
When we determine the first basis $V_1=\orth(Z_1(g_0))$, we solve a second Lyapunov equation in an arbitrary parameter $g_0^{\rr}\in\cD$ with $g_0^{\rr}\neq g_0$ to obtain the solution $Z_1(g_0^{\rr})$. 
To limit the possibilities of choosing $g_0^{\rr}$, we again define a finite subset $\cD_{\Test}\subset\cD$ and pick the parameter $g_0^{\rr}$ from this finite set $\cD_{\Test}$.
However, we  can select the parameter $g^{\rr}_0$ and the following parameters $g^{\rr}$ randomly if we want to avoid confinement to a parameter set $\cD$.
We use the solution $Z_1(g_0^{\rr})$ and the basis $V_1$ to obtain the first error equation basis 
\[
V_{1,\err}=\orth([V_1,~Z_1(g_0^{\rr})])=\orth([Z_1(g_0),~Z_1(g_0^{\rr})]).
\]
With the bases $V_1$ and $V_{1,\err}$, we define the error estimate $\Delta(g)$ as in \eqref{eq:ErrEst}.
The basis computation is described in Steps \ref{S1} to \ref{S5} of \Cref{algo:adaptiveRBM}.
After computing the first bases $V_1$ and $V_{1,\err}$, we start the optimization process of the reduced problem defined by \Cref{algo:hybridRedOpt}.
This optimization yields either the minimizer $g^*$ or, if $\mathsf{conv} = \mathsf{false}$ holds, the information that the optimization process did not converge and we need to enrich the bases.
If the bases need to be expanded, in Steps \ref{S8} to \ref{S12} 
we enlarge the bases $V_1$ and $V_{1,\err}$ as 
\[
V_1 = \orth([V_1,~Z_1(g) ]), \;\text{  and  }\;V_{1,\err}=\orth([V_{1,\err},~Z_1(g),~Z_1(g^{\rr})]),
\]
where we choose in Step \ref{S12} the parameter $g^{\rr}\in\cD_{\Test}$ that results in the largest residual 
\[
g^{\rr} = \argmax_{g\in\cD_{\Test}} \|\cR_{\rr}(g)\|_{\Fr}
\]
with $\cR_{\rr}(g)$ defined as in \eqref{eq:Resisual}.
Afterwards, in Step \ref{S13} we compute the approximated energy response value and proceed with the minimization process.

\begin{algorithm}
	\caption{Reduced energy response}\label{algo:hybridRedOpt}
	\textbf{Input:} $\cA:\cD\to\Rnn$ asymptotically stable, $\cB\in\Rnm$, $g\in\cD$, basis $V_1$, tolerance $\tol_f$.\\
	\textbf{Output:} Energy response $\tilde{\bJ}(g)$, variable $\mathsf{conv}$ that shows whether the algorithm converged.
	\begin{spacing}{1.2}	
		\begin{algorithmic}[1] 
			\State Set $\mathsf{conv} = \mathsf{true}$.
			\If{$\Delta(g)>\tol_f$}\label{S2}
			\State Set $\mathsf{conv} = \mathsf{false}$, $\tilde{\bJ}(g) = \infty$.
			\Else
			\State Solve the reduced Lyapunov equation \eqref{eq:redLE} to obtain $\widehat{P}_{11}(g)$.\label{Sf1}
			\State Set $\tilde{\bJ}(g) = \mathrm{tr}(\cC V_1\widehat{P}_{11}(g)V_1^{\T}\cC^{\T})$.\label{Sf2}
			\EndIf
		\end{algorithmic}
	\end{spacing}
\end{algorithm}

\begin{algorithm}
	\caption{Adaptive reduced basis method}\label{algo:adaptiveRBM}
	\textbf{Input:} $\cA:\cD\to\Rnn$ asymptotically stable, $\cB\in\Rnm$, initial parameters $g,\,g^\rr\in\cD_{\Test}$ with $g\neq g^{\rr}$, tolerance $\tol_f$.\\
	\textbf{Output:} Minimizer $g^{\mathrm{opt}}$, energy response $\tilde{\bJ}(g^{\mathrm{opt}})$.
	\begin{spacing}{1.2}	
		\begin{algorithmic}[1] 
			\State Choose $g_0,~ g_0^{\rr}\in\cD_{\Test}$, $g_0 \neq g_0^{\rr}$.\label{S1}
			\State Solve the Lyapunov equation from \eqref{eq:Lyap} in $g_{0}$ to obtain $Z_1(g_{0})$.
			\State Set $V_1 := \mathrm{orth}(Z_1(g_{0}))$.
			\State Solve the Lyapunov equation from \eqref{eq:Lyap} in $g_0^{\rr}$ to obtain $Z_1(g_0^{\rr})$.
			\State Set $V_{1,\err} = \orth([Z_1(g_{0}), ~Z_1(g_0^{\rr})])$.\label{S5}
			\State Apply an optimization method to optimize the function given by \Cref{algo:hybridRedOpt} to obtain the minimizer $g^{\mathrm{opt}}$, $\tilde{\bJ}(g^{\mathrm{opt}})$ and $\mathsf{conv}$.\label{S6}
			\While{$\mathsf{conv} = \mathsf{false}$}
			\State Solve the Lyapunov equation from \eqref{eq:Lyap} in $g^{\mathrm{opt}}$ to obtain $Z(g^{\mathrm{opt}})$.\label{S8}
			\State Set $V_1 := \orth([V_1, Z_1(g^{\mathrm{opt}})])$.
			\State Determine $g^{\rr}:= \argmax_{g\in\cD_{\Test}}\|R_{\rr}(g)\|_{\Fr}$.\label{S10}
			\State Solve the Lyapunov equation from \eqref{eq:Lyap} in $g^{\rr}$ to obtain $Z_1(g^{\rr})$.
			\State Set $V_{1,\err} = \orth([V_{1,\err},~Z_1(g^{\mathrm{opt}}),~ Z_1(g^{\rr})])$.\label{S12}
			\State Apply an optimization method to optimize the function given by \Cref{algo:hybridRedOpt} to obtain $g^{\mathrm{opt}}$, $\tilde{\bJ}(g^{\mathrm{opt}})$ and $\mathsf{conv}$.\label{S13}
			\EndWhile
		\end{algorithmic}
	\end{spacing}
\end{algorithm}

\begin{remark}
As in the previous section, we solve the Lyapunov equation \eqref{eq:Lyap} in $g=0\in\R^{\ell}$ (undamped system) to obtain $Z_1(0)$.
The vectors of $Z_1(0)$ are then added to the basis $V_1$ which turns out to lead to a more robust basis.
\end{remark}


\section{Implementation details}\label{sec:ImplDetail}
In this section, we specify some implementational details. 
First, in \Cref{ssec:ModRep}, we describe a transformation that leads to a numerically advantageous system. Then, in \Cref{ssec:Sign} we make use of this structure that accelerates the sign-function method that is used to solve the Lyapunov equations.

\subsection{Modal representation}\label{ssec:ModRep}
In order to simplify the computations and the numerical effort, we describe briefly an useful transformation, that is used in the following.
As shown 
in \cite{TruV07}, there exists a transformation $\Phi$, called modal matrix, such that 
\begin{align*}
\Phi^\T M\Phi = I,\qquad \Phi^\T K\Phi = \Omega^2 = \diag{\omega_1^2,\dots, \omega_n^2}.
\end{align*}
The values $\omega_1, \dots, \omega_{n}$ are the eigenvalues of the undamped system and are called \emph{eigenfrequencies}, see \cite{morTruV09}.
The transformation matrix $\Phi$ is given by the spectral decomposition of $M^{-\frac{1}{2}}KM^{-\frac{1}{2}}$ as
\begin{align*}
\Phi:= M^{-\frac{1}{2}}U, \qquad U\Omega U^\T = M^{-\frac{1}{2}}KM^{-\frac{1}{2}},
\end{align*}
where $\Omega$ is as above. 
It holds that 
$
\Phi^\T D_\intern\Phi = 2\alpha \Omega.
$
That means that $\Phi$ diagonalizes the internal damping $D_\intern$. 
Hence, this damping is called \emph{modal damping}.
The transformed mass matrix is the identity matrix, the transformed stiffness and internal damping matrix are diagonal matrices and the external damping matrix is written using its low rank factors as $(\Phi^\T  F) G(g)(\Phi^\T  F)^\T$.
Hence, the second order system \eqref{eq:SOsys} is equivalent to the first order system  \eqref{eq:FOsys} with the matrices 
\begin{align}\label{eq:decA}
\cA(g):= \widetilde{\cA} - \cU G(g)\cU^\T, \quad \widetilde{\cA} = \begin{bmatrix}
0 & I \\
-\Omega^2& -2\alpha\Omega
\end{bmatrix},~
\cU = \begin{bmatrix}
0\\
\Phi^\T  F
\end{bmatrix}, \quad \cB := \begin{bmatrix}
0\\
\Phi^\T B
\end{bmatrix}, \quad \cC := \begin{bmatrix}
C\Phi & 0
\end{bmatrix}.
\end{align}

\subsection{Solving Lyapunov equations using sign function method}\label{ssec:Sign}
As described in the previous sections, we need to solve Lyapunov equations from \eqref{eq:Lyap} in order to compute the desired energy response from \eqref{eq:SysResTr} for specific parameters $g\in\cD$.
Hence, in this section, we briefly review some numerical methods to solve the parameter-independent Lyapunov equation 
\begin{equation}\label{eq:Lyap_paraIndep}
\cA \cP + \cP \cA^{\T} = -\cB\cB^{\T}.
\end{equation}
There are multiple methods for solving this kind of equations.
If the matrix dimensions are sufficiently small, dense direct solvers such as the Hammarling method \cite{Ham82b} or the Bartels-Steward algorithm \cite{BarS72} are available. There are also dense iterative solvers such as the sign function method introduced in \cite{BenQ99}.
However, these methods are unfeasible if the matrix dimensions are large.
In this case, the alternating-direction implicit (ADI) method \cite{Pen00b, Kue16} and Krylov subspace methods \cite{SimD09} are the state of the art.
Moreover, for certain matrix structures, the sign function method can be applied, which exploits these structures to solve the Lyapunov equations efficiently.
In this work, the latter method is the one of choice, since numerical experiments indicate that it leads to the fastest results in our setting.

Within the sign-function method, which was derived in \cite{morDen19}, the structure presented in \eqref{eq:decA} is used to solve the Lyapunov equations more efficiently.
The sign function method determines a low-rank factor $\cZ$, such that $\cP \approx \cZ\cZ^{\T}$ is an approximative  solution of the Lyapunov equation in \eqref{eq:Lyap_paraIndep}.
We exploit the fact, that the Lyapunov equation in \eqref{eq:Lyap_paraIndep} is equivalent to the equation
\[
\begin{bmatrix}
I & 0\\
-\cP & I
\end{bmatrix}\begin{bmatrix}
\cA^{\T} & 0\\
0 & -\cA
\end{bmatrix} \begin{bmatrix}
I & 0\\
\cP & I
\end{bmatrix} = \begin{bmatrix}
\cA^{\T} & 0\\
\cB\cB^{\T} & -\cA
\end{bmatrix}=:\cW\qquad\text{and that}\qquad \Sign(\cW)= \small
\left[\hspace*{-3pt} \begin{array}{cc}
-I & 0 \\
2\cP & I
\end{array}\hspace*{-3pt} \right]
\]
where $\Sign(\cdot)$ denotes the sign function of a matrix. 
We observe that the $\Sign$ function of $\cW$ provides the solution $\cP\approx \cZ\cZ^{\T}$ of the Lyapunov equation from \eqref{eq:Lyap_paraIndep} in its lower left block.
The $\Sign$ function can be computed by applying Newton's method:
\begin{align}\label{eq:Newton}
A_0:=A,\;\;\; A_{k+1} = \frac{1}{2}c_kA_k + \frac{1}{2c_k}A_k^{-1}\;\;\to\;\;\;\mathrm{sign}(A),
\end{align}	
where $c_k$ is an acceleration factor and can be chosen as $c_k = \sqrt{\|A_k^{-1}\|_{\Fr}\|A_k\|_{\Fr}^{-1}}.$
The convergence of this method is shown in \cite{morRob80}.

The Newton's method described in \eqref{eq:Newton} is applied to compute $\Sign(\cW)$. 
In order to improve the efficiency while computing the inverse $\cA_k^{-1}$, the decomposition presented in \eqref{eq:decA} as well as the Sherman-Morrison-Woodbury formula is applied as described in \cite{morDen19}.
The initial values are set to be
\[
\widetilde{\cA}_{0} = \widetilde{\cA}, \quad \cU_0 = \cU, \quad \cV_0 = \cU, \quad G_0 = G, \quad \cB_0 = \cB
\]
and the iteration is defined as 
\begin{align*}
\widetilde{\cA}_{k+1} &= \frac{1}{2}\left(c_k\widetilde{\cA}_k + \frac{1}{c_k} \widetilde{\cA}_k^{-1}\right),\;
\cU_{k+1} = \left[ U_k , \; \widetilde{\cA}_k^{-1}\cU_k \right],\;
\cV_{k+1} = \left[ V_k, \; \widetilde{\cA}_k^{-\T}\cV_k \right],\\
G_{k+1} &= \frac{1}{2} \mathrm{diag} \left( c_k G_k,\; -\frac{1}{c_k}(G_k^{-1} - \cV_k^{\T}\widetilde{\cA}_k^{-1}\cU_k)^{-1} \right), \; \cB_{k+1} = \frac{1}{\sqrt{2}}\begin{bmatrix}
\sqrt{c_k}\cB_{k} & \frac{1}{\sqrt{c_k}}\cA_{k}^{-1}\cB_k
\end{bmatrix}
\end{align*}
so that $\cB_{k+1}$ converges to $\frac{1}{\sqrt{2}}\cZ$ for $k\to\infty$.
We stop this method if $\|\cA_k+I\|^{2}\leq\tol$ since $\cA_k$ converges to $-I$, or if a maximum number of iterations $\mathrm{iter}_{\max}$ is exceeded.

One disadvantage of this method is the high growth rate of the dimension of the low-rank factor $\cB_{k+1}$.
Therefore, even with internal truncation techniques, the method becomes rather slow if it does not converge after a few steps. 
Hence, in this work, we set the maximum number of iterations $\mathrm{iter}_{\max}$ to be quite small, so that the method is interrupted before the low-rank factors are of too large dimensions.
Hence, the low-rank factor approximation is not guaranteed to be accurate, however, it still spans a space that is a good enough approximation of the solution space as we will show in the numerical examples.

\section{Examples and numerical results}\label{sec:Results}
In this section, we illustrate the accelerations that arise when we apply the methods presented in this paper.
Therefore, we consider two examples that are presented in \cite{morTomBG18}.
The computations have been done on a computer with 2 Intel Xeon Silver 4110 CPUs running at 2.1 GHz and equipped with 192 GB total main memory. 
The experiments use \matlab 2021a.

\subsection{Example 1}\label{sec:NumEx1}
The first example was introduced in \cite{morTomBG18} and arises in mechanical constructions with $n$ consecutive masses.
Each mass $m_j$ is connected to the direct neighbor masses $m_{j-1}$ and $m_{j+1}$ by springs with stiffness values $k_j$ and $k_{j+1}$.
Additionally, each mass is connected by springs with stiffness values $k_{j-1}$ and $k_{j+2}$ to the masses next to the neighbor masses $m_{j-2}$ and $m_{j+2}$.
The outermost masses are connected to fixed objects via springs with constants $2k_1$ and $2k_n$.
This construction results in the following mass and stiffness matrix
\begin{align*}
M &:= \diag{m_1, ~ \dots ,~m_n},\\ K &:= \begin{bmatrix}
2k_1 + 2k_2 & -k_2        & -k_3        &                     &          & \\
-k_2        & 2k_2 + 2k_3 & -k_3        & -k_4                &          & \\
-k_3        & -k_3        & 2k_3 + 2k_4 & -k_4                &  -k_5      & \\
& \ddots      & \ddots      & \ddots              &  \ddots    & \\
&             &             & 2k_{n-2} + 2k_{n-1} & -k_{n-2} & -k_{n-1} \\
&             &             & -k_{n-1}            & 2k_{n-1} + 2k_n & -k_n \\
&             &             & -k_{n-2}            & -k_n& 2k_{n} + 2k_{n+1}
\end{bmatrix}.
\end{align*}
We consider an example of dimension $n=1900$ with stiffness constants $k_j=500,~ j=1, ~\dots, ~ n$.
The mass values are chosen as 
\[
m_j = \begin{cases}
144-\frac{3}{20}j,  & j = 1, \dots, 475,\\
\frac{j}{10} + 25,  & j = 476,\dots,1900.
\end{cases}
\]
The internal damping $D_{\intern}$ is built as described in Equation \eqref{eq:intDamp} where the scaling factor is $\alpha = 0.005$.
We consider external disturbance forces that attack at the sequential masses from $m_{471}$ to $m_{480}$.
Hence, in the input matrix $B$ the values at positions $471$ to $480$ are set to be 
\[
B(471:480,\, 1:10) = \diag{
10,\, 20,\, 30,\, 40,\, 50,\, 50,\, 40,\, 30,\, 20,\, 10}.
\] The other entries of $B$ are equal to zero.
Consequently, we have an $(n\times 10)$-dimensional input matrix $B$ where the highest magnitude of disturbance is applied to the mass in the center, whereby  the disturbance magnitude gets smaller in the outer masses.
To observe the system behavior, we consider the displacement of the states $x_{100}(t),~ x_{200}(t),~\dots, x_{1800}(t)$.
Hence, the output matrix $C$ is $18\times n$-dimensional and has zero entries everywhere except at the positions $(1, 100),~(2, 200),~\dots,~(18, 1800)$ where the entries are equal to one.

Now, we consider the external dampers that we want to optimize.
We consider four dampers at the positions $j,~ j+1,~k,~k+1$ where $j$ and $k$ can take the following values
\[
\{ (j, k)~ |~j \in \{50,~ 150,~ 250,~ 350\},~ k\in\{850,~950,~\dots,~1850\}  \}
\]
such that we obtain $44$ possible damping configurations.
For each damping configuration, we optimize the damping values individually.
The damping gains $g$ consist of two values $g_1$ and $g_2$, where the dampers at the $j$-th and the $(j+1)$-st position have the damping value $g_1$ and the dampers at the $k$-th and the $(k+1)$-st position the damping value $g_2$.
We assume that the damping values $g_1$ and $g_2$ lie in the interval $[500, 4000]$.

To optimize the damping gains for the different damping configurations, we use the \matlab-function \texttt{fminfun}, where we stop the minimization process if the difference between two successive function values or damping gains is smaller than a tolerance of $10^{-4}$. We start the optimization process at $g_0 = \begin{bmatrix}
1000 & 1000
\end{bmatrix}^{\T}$ for all damping configurations.
To solve the Lyapunov equations from \eqref{eq:Lyap}, we use the sign-function method that is presented in \Cref{ssec:Sign} with $\tol = 10^{-6}$ and a maximum iteration number of $\mathrm{iter}_{\max}=10$ because of the fast dimension growth within the method.
The tolerance for the function value error that indicates whether a basis is sufficiently detailed is set to be $\tol_{f}=10^{-3}$.
As test-parameter set $\cD_{\Test}$ for the RBM, we use 36 uniformly distributed parameters in $[500, 4000]\times [500, 4000]$.

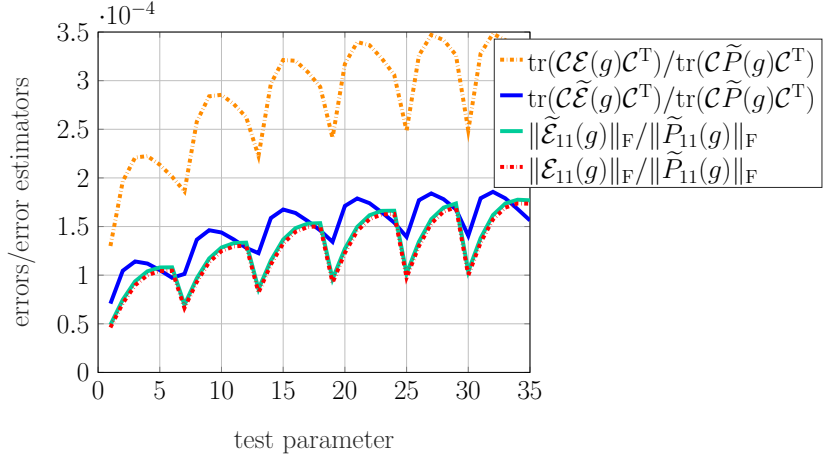
\begin{figure}[h!]
	\centering
%
%
\definecolor{mycolor1}{rgb}{1.0, 0.55, 0.0}%
\definecolor{mycolor2}{rgb}{0.0, 0.8, 0.6}%
\begin{tikzpicture}[scale=0.5]
\huge
\begin{axis}[%
width=4.521in,
height=3.566in,
at={(0.758in,0.481in)},
scale only axis,
xmin=0,
xmax=35,
xlabel style={font=\color{white!15!black},yshift=-0.8cm},
xlabel={test parameter},
ymin=0,
ymax=0.00035,
yminorticks=true,
ylabel style={font=\color{white!15!black},yshift=0.8cm},
ylabel={errors/error estimators},
axis background/.style={fill=white},
legend style={legend cell align=left, align=left, draw=white!15!black,xshift=8cm},
ymajorgrids,
yminorgrids,
xmajorgrids,
]
\addplot [color=mycolor1, dashdotted, line width=3.0pt]
  table[row sep=crcr]{%
1	0.000129979149544616\\
2	0.000196447945154736\\
3	0.000221007517890273\\
4	0.000222414833480043\\
5	0.000213540417391122\\
6	0.000201166842684884\\
7	0.000186102692577881\\
8	0.000257857826667656\\
9	0.000283923465664816\\
10	0.000285431846644738\\
11	0.000275903948381518\\
12	0.000262471553478238\\
13	0.000222498874327057\\
14	0.000296740038076542\\
15	0.000321242526370298\\
16	0.000320435995185731\\
17	0.000308803433547239\\
18	0.000293659324594416\\
19	0.000241126671427496\\
20	0.000317110758832874\\
21	0.000339709952382304\\
22	0.000336417485210556\\
23	0.000322666261915051\\
24	0.000305912966941398\\
25	0.000248275596899045\\
26	0.000325859682150337\\
27	0.000347062083596779\\
28	0.000341800915178897\\
29	0.000326349924287322\\
30	0.00024888285155549\\
31	0.000327936865284042\\
32	0.000348296378057114\\
33	0.000341636769914862\\
34	0.000324935321620265\\
35	0.000305942321620584\\
};
\addlegendentry{$\mathrm{tr}(\mathcal{CE}(g)\mathcal{C}^{\mathrm{T}})/\mathrm{tr}(\mathcal{C}\widetilde{P}(g)\mathcal{C}^{\mathrm{T}})$}

\addplot [color=blue, line width=3.0pt]
  table[row sep=crcr]{%
1	7.07434692456564e-05\\
2	0.000104560196301536\\
3	0.000114096736675388\\
4	0.000111935076549906\\
5	0.000105187827330487\\
6	9.72832528730167e-05\\
7	0.000101233410199501\\
8	0.000136483052755977\\
9	0.000146212453116764\\
10	0.000143883515515387\\
11	0.000136669354196733\\
12	0.000128099409827194\\
13	0.00012259833260378\\
14	0.000158726606210735\\
15	0.000167508258626532\\
16	0.000163932520890155\\
17	0.000155564589878687\\
18	0.00014601548594629\\
19	0.000134257956383023\\
20	0.00017113783857969\\
21	0.000178871939908076\\
22	0.000173978404637978\\
23	0.000164471841496892\\
24	0.000154033322468292\\
25	0.000139338945617286\\
26	0.000177054311276035\\
27	0.000184048870227697\\
28	0.000178115532862901\\
29	0.000167699800713842\\
30	0.000140536386088681\\
31	0.00017911459199608\\
32	0.000185688765898028\\
33	0.000179020382491986\\
34	0.000167936716108263\\
35	0.000156280238306055\\
};
\addlegendentry{$\mathrm{tr}(\mathcal{C\widetilde{E}}(g)\mathcal{C}^{\mathrm{T}})/\mathrm{tr}(\mathcal{C}\widetilde{P}(g)\mathcal{C}^{\mathrm{T}})$}

\addplot [color=mycolor2, line width=3.0pt]
  table[row sep=crcr]{%
1	4.89744483731591e-05\\
2	7.47190554933174e-05\\
3	9.36483853923213e-05\\
4	0.000103965845704256\\
5	0.00010792155394643\\
6	0.000108051057863666\\
7	6.93863255275443e-05\\
8	9.6778453603665e-05\\
9	0.000117161868457667\\
10	0.000128421367713344\\
11	0.000132968461686857\\
12	0.000133456249096494\\
13	8.56190168153358e-05\\
14	0.000114963633579455\\
15	0.000136588292469912\\
16	0.000148367037906332\\
17	0.000153037203824167\\
18	0.000153473816049571\\
19	9.60187797452704e-05\\
20	0.000127152435044403\\
21	0.000149657666444775\\
22	0.000161610451172732\\
23	0.000166111202313845\\
24	0.000166266951746393\\
25	0.000101617714394436\\
26	0.000134355518980116\\
27	0.000157516730371608\\
28	0.000169475645747829\\
29	0.000173681281217875\\
30	0.000103951302840006\\
31	0.00013809046100624\\
32	0.000161774896676832\\
33	0.000173685326892759\\
34	0.000177577765377314\\
35	0.000176972875997125\\
};
\addlegendentry{$\|\widetilde{\mathcal{E}}_{11}(g)\|_{\mathrm{F}}/\|\widetilde{P}_{11}(g)\|_{\mathrm{F}}$}

\addplot [color=red, dashdotted, line width=3.0pt]
table[row sep=crcr]{%
1	4.6223604473297e-05\\
2	7.09677593704801e-05\\
3	8.96087577607888e-05\\
4	0.00010001104173973\\
5	0.000104216325179663\\
6	0.000104641869437528\\
7	6.61979272554977e-05\\
8	9.30096558497201e-05\\
9	0.000113326810938973\\
10	0.000124762273417749\\
11	0.000129581210384118\\
12	0.000130355518131022\\
13	8.19947997243143e-05\\
14	0.000110911079001516\\
15	0.000132571699507988\\
16	0.000144587110727574\\
17	0.000149563818503725\\
18	0.000150306290392719\\
19	9.20721764003681e-05\\
20	0.000122818410511462\\
21	0.000145407783279944\\
22	0.00015763675470542\\
23	0.000162473954916246\\
24	0.000162957390340622\\
25	9.74840381753909e-05\\
26	0.000129824813611788\\
27	0.000153089659003888\\
28	0.000165349420370901\\
29	0.000169913899662142\\
30	9.97366946554048e-05\\
31	0.00013344616950576\\
32	0.000157237334411488\\
33	0.000169462846665424\\
34	0.000173729784856932\\
35	0.000173484202379791\\
};
\addlegendentry{$\|\mathcal{E}_{11}(g)\|_{\mathrm{F}}/\|\widetilde{P}_{11}(g)\|_{\mathrm{F}}$}
\end{axis}

\end{tikzpicture}%
	\caption{Errors estimator for Example 1 for the first damping configuration and the first step of the RBM.}
	\label{pic:ErrorsEst1}
\end{figure}

After the first step of the RBM for the $11$-th damper configuration, we display the quality of the error estimator in \Cref{pic:ErrorsEst1}.
We observe that the relative error in the position controllability Gramian $\|\cE_{11}(g)\|_{\mathrm{F}}/\|\widetilde{P}_{11}(g)\|_{\mathrm{F}}$ and the corresponding estimate $\|\widetilde{\cE}_{11}(g)\|_{\mathrm{F}}/\|\widetilde{P}_{11}(g)\|_{\mathrm{F}}$ are very close, so the error is well approximated.
On the other hand, the energy response is underestimated since we are considering an error estimator and not an error bound.
However, for our purpose the quality was good enough, since the energy response error and its approximate are of a similar magnitude.

\begin{figure}[h!]
	\centering
%
%
\definecolor{mycolor1}{rgb}{0.85098,0.32549,0.09804}%
\definecolor{mycolor2}{rgb}{0.00000,0.49804,0.00000}%
\begin{tikzpicture}[scale=0.4]
\huge
\begin{axis}[%
width=9.155in,
height=2.422in,
at={(1.536in,0.584in)},
scale only axis,
xmin=0,
xmax=45,
xlabel style={font=\color{white!15!black},yshift=-0.8cm},
xlabel={Configuration},
xtick={0, 4,8,12,16,20,24,28,32,36,40,44},
ymode=log,
ymin=0.00001,
ymax=0.01,
yminorticks=true,
ylabel style={font=\color{white!15!black},yshift=0.8cm},
ylabel={Error},
axis background/.style={fill=white},
ymajorgrids,
yminorgrids,
legend style={legend cell align=left, align=left, draw=white!15!black,xshift=4.8cm}
]
\addplot [color=mycolor1, line width=3.0pt]
  table[row sep=crcr]{%
1	0.000626593735072787\\
2	0.000694325640349901\\
3	0.000417091910801606\\
4	0.000393048753844531\\
5	0.000345615178982798\\
6	0.000126755989779278\\
7	0.000616984558241266\\
8	0.000347891051556885\\
9	0.000664264667149432\\
10	0.00069676811590851\\
11	8.35253657013123e-05\\
12	0.00294337238795238\\
13	0.000453546870230545\\
14	0.000378852017708773\\
15	0.000490270530286499\\
16	0.000406920654682031\\
17	0.000734755697417131\\
18	0.000413870603886496\\
19	0.00101379153035266\\
20	0.00129890948801632\\
21	0.000177307396859872\\
22	0.000210818834412987\\
23	0.000954377170274269\\
24	0.00101355767152458\\
25	0.000753273407085916\\
26	0.00106710382155261\\
27	0.000610205432098413\\
28	0.000635589096666949\\
29	8.78020428990235e-05\\
30	0.000698463520549012\\
31	0.000824727125516442\\
32	0.000682895624905053\\
33	0.000172435278440918\\
34	0.00579867251862152\\
35	0.00191021218593527\\
36	0.00394081640648658\\
37	0.000970980472135939\\
38	0.000975179080286559\\
39	0.00290143670990187\\
40	0.000273462730469428\\
41	0.000912927069667753\\
42	0.000678635698622112\\
43	0.00112771539738795\\
44	0.000116407849319294\\
};
\addlegendentry{Adaptive RBM}

\addplot [color=blue, dashdotted, line width=3.0pt]
  table[row sep=crcr]{%
1	0.000626593735072787\\
2	0.000694325640349901\\
3	0.000417091910801606\\
4	0.000393048753844531\\
5	0.000345615178982798\\
6	0.000126755989779278\\
7	0.000616984558241266\\
8	0.000347891051556885\\
9	0.000664264667149432\\
10	0.00069676811590851\\
11	8.35253657013123e-05\\
12	0.00294337238795238\\
13	0.000453546870230545\\
14	0.000378852017708773\\
15	0.000490270530286499\\
16	0.000406920654682031\\
17	0.000734755697417131\\
18	0.000413870603886496\\
19	0.00101379153035266\\
20	0.00129890948801632\\
21	0.000177307396859872\\
22	0.000210818834412987\\
23	0.000954377170274269\\
24	0.00101355767152458\\
25	0.000753273407085916\\
26	0.00106710382155261\\
27	0.000610205432098413\\
28	0.000635589096666949\\
29	8.78020428990235e-05\\
30	0.000698463520549012\\
31	0.000824727125516442\\
32	0.000682895624905053\\
33	0.000172435278440918\\
34	0.00579867251862152\\
35	0.00191021218593527\\
36	0.00394081640648658\\
37	0.000970980472135939\\
38	0.000975179080286559\\
39	0.00290143670990187\\
40	0.000273462730469428\\
41	0.000912927069667753\\
42	0.000678635698622112\\
43	0.00112771539738795\\
44	0.000116407849319294\\
};
\addlegendentry{RBM}

\addplot [color=mycolor2, dotted, line width=3.0pt]
  table[row sep=crcr]{%
1	0.00051344003425912\\
2	6.24966846341575e-05\\
3	0.000387389410002308\\
4	0.000253733404817875\\
5	0.000521199206469401\\
6	0.000112306526781427\\
7	0.000885210231417762\\
8	0.000113957522137396\\
9	0.000163291968651406\\
10	6.90734589138372e-06\\
11	0.00012027632572108\\
12	0.000828775104540662\\
13	0.000517841406367612\\
14	0.000132648117995703\\
15	0.000453517726023542\\
16	0.000601792311426158\\
17	0.000585138277213544\\
18	0.000561842499251187\\
19	0.000615904694115323\\
20	5.83976869889715e-05\\
21	0.000426811994526225\\
22	0.000275033617280864\\
23	0.000152848182321302\\
24	0.000463668494033791\\
25	0.000249962397101177\\
26	0.000303403863163307\\
27	0.000612581894994207\\
28	0.000679547994014148\\
29	0.000228657217789434\\
30	0.000649540627240008\\
31	0.000145454748727495\\
32	0.000635348344211983\\
33	0.000212373316967076\\
34	0.00295914908703432\\
35	0.000145396373811522\\
36	0.000498502961387953\\
37	0.000508306494796029\\
38	0.000970552790465568\\
39	0.000237110897159722\\
40	0.000611814127549196\\
41	0.000148696014179486\\
42	0.000139084223266099\\
43	0.000200741922059587\\
44	0.000183322290027637\\
};
\addlegendentry{sym2IRKA}

\end{axis}

\begin{axis}[%
width=5.833in,
height=4.375in,
at={(0in,0in)},
scale only axis,
xmin=0,
xmax=1,
ymin=0,
ymax=1,
axis line style={draw=none},
ticks=none,
axis x line*=bottom,
axis y line*=left,
legend style={legend cell align=left, align=left, draw=white!15!black}
]
\end{axis}
\end{tikzpicture}%
	\caption{Errors example 1}
	\label{pic:ErrorsEx1}
\end{figure}
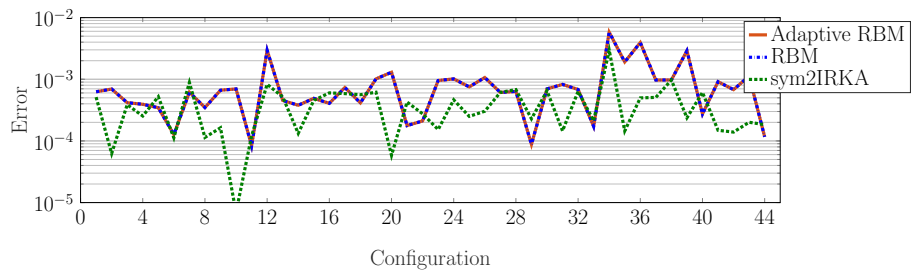

Since the initial value $g_0$ is known, we choose this parameter as first one that is evaluated within the RBM.
The first parameter $g_0^{\rr}$ that is used to obtain a first error equation basis is chosen to be $g_0^{\rr} = \begin{bmatrix}
100 & 100
\end{bmatrix}$ within the RBM and the adaptive RBM.
Within the state of the art methods, the  symmetric IRKA approach for second-order systems (\texttt{sym2IRKA}) leads to the highest acceleration rates in the optimization of external dampers and their viscosities. Hence, we will compare in the following our approach to the \texttt{sym2IRKA} method from \cite{morTomBG18}.
Since its code is not available, we implemented the method ourselves and compare the results to the best implementation and configurations we were able to generate using the \texttt{sym2IRKA} method. 
This method is a projection method as well that enriches the corresponding basis by vectors that are generated using an IRKA approach.
This method starts the optimization processes at $g_0$ as well. 
For every enriched basis a new optimization process is started.
In every iteration, $60$ vectors are added to the basis and the method stops if the relative error between two consecutive optima is smaller than the tolerance $\tol_{\mathrm{IRKA}}=10^{-3}$ .
The relative errors between the optimal damping gain and the approximations obtained using the methods presented in the previous sections are presented in \Cref{pic:ErrorsEx1}.
We observe that all errors are smaller than $10^{-2}$ and hence that the damping values are for our purposes sufficiently detailed.

\begin{figure}[h!]
	\centering
%
%
\definecolor{mycolor1}{rgb}{0.85098,0.32549,0.09804}%
\definecolor{mycolor2}{rgb}{0.00000,0.49804,0.00000}%
\begin{tikzpicture}[scale=0.4]
\huge

\begin{axis}[%
width=9.155in,
height=2.422in,
at={(1.536in,0.584in)},
scale only axis,
xmin=0,
xmax=45,
xlabel style={font=\color{white!15!black},yshift=-0.8cm},
xlabel={Configuration},
xtick={0, 4,8,12,16,20,24,28,32,36,40,44},
ymin=0,
ymax=500,
ylabel style={font=\color{white!15!black},yshift=0.8cm},
ylabel={Acceleration rates},
axis background/.style={fill=white},
ymajorgrids,
legend style={legend cell align=left, align=left, draw=white!15!black,xshift=4.8cm}
]
\addplot [color=blue, dashdotted, line width=3.0pt]
  table[row sep=crcr]{%
1	240.909108086108\\
2	325.650731113911\\
3	327.661066362465\\
4	295.969821573528\\
5	305.564183351121\\
6	309.687887459196\\
7	391.318657379143\\
8	360.032551521148\\
9	423.446799504427\\
10	409.774549420301\\
11	465.286505248661\\
12	473.185522641623\\
13	266.979068116091\\
14	255.440221798512\\
15	343.504008665095\\
16	314.30562960231\\
17	337.051932877097\\
18	332.170911175488\\
19	340.056759560642\\
20	362.539694054164\\
21	367.430027074702\\
22	368.618869413646\\
23	294.904041937161\\
24	399.130074441706\\
25	273.041875241859\\
26	274.061525787127\\
27	285.222071699998\\
28	330.674075700646\\
29	321.072152143809\\
30	313.057278744661\\
31	329.563513439663\\
32	389.888906395717\\
33	316.315259218701\\
34	451.050615148744\\
35	314.652812471661\\
36	304.856453331912\\
37	273.222015850225\\
38	254.935276611861\\
39	373.386479354171\\
40	259.534539732825\\
41	353.726947603028\\
42	373.937776894021\\
43	421.276236977842\\
44	360.731196872534\\
};
\addlegendentry{RBM}

\addplot [color=mycolor1, line width=3.0pt]
  table[row sep=crcr]{%
1	175.503248840924\\
2	201.539291391664\\
3	197.513602615469\\
4	190.58856681484\\
5	200.042526171592\\
6	194.983539779988\\
7	245.411121720024\\
8	217.251270662793\\
9	255.905900532808\\
10	240.79954748297\\
11	266.237063431131\\
12	279.046361998329\\
13	167.72234592823\\
14	172.469065830092\\
15	193.506444368773\\
16	192.666999848025\\
17	211.734440378211\\
18	203.480091956605\\
19	206.14914173581\\
20	219.172040229368\\
21	225.66664130203\\
22	215.357539229028\\
23	191.606989520865\\
24	242.436099058231\\
25	172.863926218824\\
26	177.773445203042\\
27	179.661999868256\\
28	216.248597646515\\
29	202.929598814485\\
30	192.227454374358\\
31	201.189229024314\\
32	236.914498600259\\
33	196.402125878265\\
34	247.26273813217\\
35	195.579663330877\\
36	188.226255908431\\
37	166.952874481076\\
38	156.358388532908\\
39	211.435986647993\\
40	164.514695540912\\
41	220.476165233916\\
42	227.095594978097\\
43	217.711086416188\\
44	215.454260858486\\
};
\addlegendentry{Adaptive RBM}

\addplot [color=mycolor2, dotted, line width=3.0pt]
  table[row sep=crcr]{%
1	64.4561897264377\\
2	154.917853717833\\
3	136.138680097926\\
4	135.376178949809\\
5	72.9701507710552\\
6	100.111082980348\\
7	136.244982693742\\
8	93.3695909699896\\
9	125.058131020258\\
10	212.418415404253\\
11	168.300262574071\\
12	19.6209559744924\\
13	70.5619172131191\\
14	113.215216234155\\
15	20.7023105903495\\
16	51.0658560384228\\
17	174.974764355487\\
18	124.708626494749\\
19	160.128715071969\\
20	187.623525659378\\
21	118.567702634769\\
22	187.257665472557\\
23	27.3455665771295\\
24	62.7749214541986\\
25	26.4350063411814\\
26	50.140135317725\\
27	80.431223826172\\
28	133.216667228328\\
29	143.05948583407\\
30	123.051134775172\\
31	70.3542313057876\\
32	125.266304051127\\
33	130.712993293408\\
34	173.980114057221\\
35	71.6986751716813\\
36	159.621652198273\\
37	34.9784440575198\\
38	92.0249016937064\\
39	93.3256481094097\\
40	89.0606104609583\\
41	158.173178307632\\
42	139.274769955606\\
43	139.286241743159\\
44	134.860812729911\\
};
\addlegendentry{sym2IRKA}

\end{axis}

\end{tikzpicture}%
	\caption{Acceleration rates example 1}
	\label{pic:AccelerationsEx1}
\end{figure}
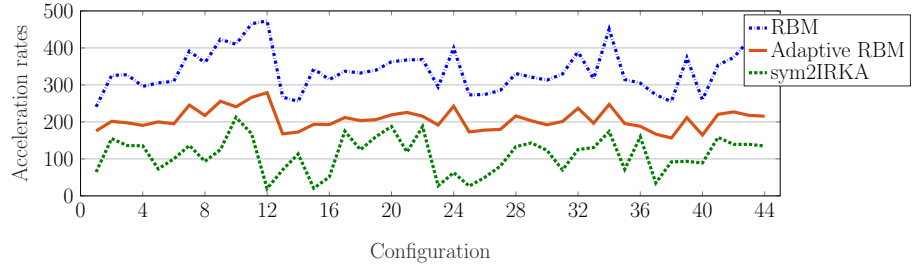

Now, we evaluate the optimization times, that include for the RBM the offline and the online phase. 
Offline, we determine a low-rank factor of the solution of the Lyapunov equation \eqref{eq:Lyap} for the undamped system. 
This low-rank factor is included into the bases in the RBM, the adaptive RBM, and the \texttt{sym2IRKA} approach.
Since this low-rank factor is computed beforehand, the computation time of $5.2$ seconds is not taken into account in any of these methods.
We compare the times in \Cref{pic:TimesEx1} and the acceleration rates for the different methods in \Cref{pic:AccelerationsEx1}.
The \matlab-solver \texttt{lyapchol} is used to solve the Lyapunov equations from \eqref{eq:Lyap}.
We observe, that the reduced optimization procedures lead to significantly faster results so that the acceleration rates are in average $339$ for the reduced basis method and $208$ for the adaptive procedure.
These results are comparable with those from \cite{morTomBG18}.
Here the acceleration rate is in average $78$ in our \texttt{sym2IRKA} implementation, that is less improvement than with our method.
However, we note that in \cite{morTomBG18} the acceleration rate for this example was $346$.
We also observe that the adaptive reduced basis method is slower than the one where offline and online phase are decoupled. 
This is because already the first basis that is equal for both methods is sufficiently good while for the adaptive method there are additional computational cost for the evaluations of the error estimators.
We want to mention, that the adaptive method can still be advantageous since we do not need a parameter set $\cD$ in advanced for it.
Only for the error estimator the test-parameter set $\cD_{\Test}\subset\cD$ is needed that can be replaced by choosing arbitrary parameters in a certain range around the damping values that are attained during the optimization process.

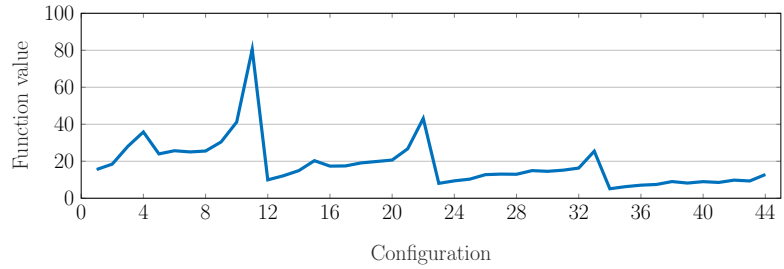
\begin{figure}[h!]
	\centering
%
%
\definecolor{mycolor1}{rgb}{0.00000,0.44700,0.74100}%
\definecolor{mycolor2}{rgb}{0.00000,0.44706,0.74118}%
\begin{tikzpicture}[scale=0.4]
\huge

\begin{axis}[%
width=9.155in,
height=2.422in,
at={(1.536in,0.584in)},
scale only axis,
xmin=0,
xmax=45,
xlabel style={font=\color{white!15!black},yshift=-0.8cm},
xlabel={Configuration},
xtick={0, 4,8,12,16,20,24,28,32,36,40,44},
ymin=0,
ymax=100,
ylabel style={font=\color{white!15!black},yshift=0.8cm},
ylabel={Function value},
axis background/.style={fill=white},
ymajorgrids
]
\addplot [color=mycolor1, forget plot]
  table[row sep=crcr]{%
1	15.5589587970499\\
2	18.4717799011823\\
3	28.0528329045287\\
4	35.8319360964957\\
5	23.9811463064594\\
6	25.6909959897258\\
7	25.0446900592919\\
8	25.5307787312104\\
9	30.4472705608641\\
10	41.1793856075722\\
11	80.4622662213073\\
12	9.94796647068045\\
13	12.1355816064576\\
14	14.9460899850742\\
15	20.323053426102\\
16	17.3685769710211\\
17	17.4736601036392\\
18	19.0995795027775\\
19	19.8724741818882\\
20	20.6522496972096\\
21	26.7318639673428\\
22	43.1245220640656\\
23	8.04986916197914\\
24	9.39909622772578\\
25	10.3430685756313\\
26	12.7594888199609\\
27	13.0549363293556\\
28	12.9681722847333\\
29	14.9372472076993\\
30	14.5631909976193\\
31	15.1755115463375\\
32	16.3202140296766\\
33	25.436453006106\\
34	5.14116840631789\\
35	6.26117463835468\\
36	7.07183925291042\\
37	7.45607440549751\\
38	9.01524276733949\\
39	8.21098598775817\\
40	8.97112315635505\\
41	8.54728441201632\\
42	9.78912639305556\\
43	9.34101652534681\\
44	12.8105301088211\\
};
\addplot [color=mycolor2, line width=3.0pt, forget plot]
  table[row sep=crcr]{%
1	15.5589587970499\\
2	18.4717799011823\\
3	28.0528329045287\\
4	35.8319360964957\\
5	23.9811463064594\\
6	25.6909959897258\\
7	25.0446900592919\\
8	25.5307787312104\\
9	30.4472705608641\\
10	41.1793856075722\\
11	80.4622662213073\\
12	9.94796647068045\\
13	12.1355816064576\\
14	14.9460899850742\\
15	20.323053426102\\
16	17.3685769710211\\
17	17.4736601036392\\
18	19.0995795027775\\
19	19.8724741818882\\
20	20.6522496972096\\
21	26.7318639673428\\
22	43.1245220640656\\
23	8.04986916197914\\
24	9.39909622772578\\
25	10.3430685756313\\
26	12.7594888199609\\
27	13.0549363293556\\
28	12.9681722847333\\
29	14.9372472076993\\
30	14.5631909976193\\
31	15.1755115463375\\
32	16.3202140296766\\
33	25.436453006106\\
34	5.14116840631789\\
35	6.26117463835468\\
36	7.07183925291042\\
37	7.45607440549751\\
38	9.01524276733949\\
39	8.21098598775817\\
40	8.97112315635505\\
41	8.54728441201632\\
42	9.78912639305556\\
43	9.34101652534681\\
44	12.8105301088211\\
};
\end{axis}
\end{tikzpicture}%
	\caption{Function values obtained in Example 1}
	\label{pic:FunValEx1}
\end{figure}

In \Cref{pic:FunValEx1}, the function values for all $44$ damping configurations are evaluated. We see that the optimal damping configuration is the $34$-th one that has the damping positions $j=350,~k=850$.

\begin{figure}[h!]
	\centering
%
%
\definecolor{mycolor1}{rgb}{0.85098,0.32549,0.09804}%
\definecolor{mycolor2}{rgb}{0.00000,0.49804,0.00000}%
\begin{tikzpicture}[scale=0.4]
\huge
\begin{axis}[%
width=9.155in,
height=2.422in,
at={(1.536in,0.584in)},
scale only axis,
xmin=0,
xmax=45,
xlabel style={font=\color{white!15!black},yshift=-0.8cm},
xlabel={Configuration},
xtick={0, 4,8,12,16,20,24,28,32,36,40,44},
ymode=log,
ymin=10,
ymax=100000,
yminorticks=true,
ylabel style={font=\color{white!15!black},yshift=0.8cm},
ylabel={Optimization time},
axis background/.style={fill=white},
ymajorgrids,
yminorgrids,
legend style={legend cell align=left, align=left, draw=white!15!black,xshift=4.8cm}
]
\addplot [color=black, line width=3.0pt]
  table[row sep=crcr]{%
1	11816.770988\\
2	15424.798582\\
3	15737.115726\\
4	13247.909623\\
5	14441.759911\\
6	15044.923105\\
7	18512.973017\\
8	16735.864735\\
9	20664.556547\\
10	19683.092267\\
11	22333.201818\\
12	22972.753497\\
13	11802.512128\\
14	11247.322635\\
15	16029.112705\\
16	15058.396858\\
17	16223.708969\\
18	15526.750269\\
19	15998.542229\\
20	18208.639518\\
21	16859.625631\\
22	17330.727837\\
23	13473.625117\\
24	19830.422922\\
25	12390.017763\\
26	12696.153689\\
27	13863.598996\\
28	15551.325006\\
29	15027.995273\\
30	14562.631007\\
31	15546.543781\\
32	18489.610374\\
33	14246.399597\\
34	22328.022569\\
35	14329.908946\\
36	14137.204645\\
37	12851.974011\\
38	11506.073889\\
39	18005.051125\\
40	11858.464546\\
41	16538.797736\\
42	17727.906501\\
43	21386.699442\\
44	16876.428835\\
};
\addlegendentry{Original}

\addplot [color=blue, dashdotted, line width=3.0pt]
  table[row sep=crcr]{%
1	49.050744\\
2	47.36608\\
3	48.028641\\
4	44.761015\\
5	47.262607\\
6	48.580922\\
7	47.309201\\
8	46.48431\\
9	48.800833\\
10	48.033955\\
11	47.998817\\
12	48.549147\\
13	44.207631\\
14	44.031134\\
15	46.663539\\
16	47.910045\\
17	48.134152\\
18	46.743257\\
19	47.046682\\
20	50.22523\\
21	45.885269\\
22	47.015303\\
23	45.688167\\
24	49.684111\\
25	45.37772\\
26	46.325925\\
27	48.606333\\
28	47.029163\\
29	46.805664\\
30	46.517465\\
31	47.173134\\
32	47.422766\\
33	45.03861\\
34	49.502255\\
35	45.54197\\
36	46.373316\\
37	47.038574\\
38	45.133314\\
39	48.220951\\
40	45.691277\\
41	46.755832\\
42	47.408707\\
43	50.766451\\
44	46.783946\\
};
\addlegendentry{RBM}

\addplot [color=mycolor1, line width=3.0pt]
  table[row sep=crcr]{%
1	67.330782\\
2	76.534945\\
3	79.676111\\
4	69.510516\\
5	72.193449\\
6	77.159965\\
7	75.436569\\
8	77.0346\\
9	80.750606\\
10	81.74057\\
11	83.884646\\
12	82.325938\\
13	70.369348\\
14	65.213565\\
15	82.835033\\
16	78.157634\\
17	76.62291\\
18	76.305992\\
19	77.60664\\
20	83.079208\\
21	74.710314\\
22	80.47421\\
23	70.319069\\
24	81.796494\\
25	71.674976\\
26	71.417605\\
27	77.164893\\
28	71.914108\\
29	74.055216\\
30	75.75729\\
31	77.273241\\
32	78.043389\\
33	72.536891\\
34	90.300798\\
35	73.268911\\
36	75.107506\\
37	76.979651\\
38	73.587826\\
39	85.156039\\
40	72.081491\\
41	75.013994\\
42	78.06363\\
43	98.234315\\
44	78.32952\\
};
\addlegendentry{Adaptive RBM}

\addplot [color=mycolor2, dotted, line width=3.0pt]
  table[row sep=crcr]{%
1	183.330275\\
2	99.567598\\
3	115.596212\\
4	97.859976\\
5	197.913253\\
6	150.282293\\
7	135.880035\\
8	179.243205\\
9	165.239608\\
10	92.661892\\
11	132.698556\\
12	1170.827432\\
13	167.264618\\
14	99.34462\\
15	774.266845\\
16	294.881904\\
17	92.720279\\
18	124.50422\\
19	99.910514\\
20	97.048808\\
21	142.194082\\
22	92.550165\\
23	492.716985\\
24	315.897216\\
25	468.697363\\
26	253.213391\\
27	172.365884\\
28	116.737082\\
29	105.047178\\
30	118.346174\\
31	220.975249\\
32	147.602426\\
33	108.989927\\
34	128.336636\\
35	199.862953\\
36	88.566961\\
37	367.425549\\
38	125.032178\\
39	192.927148\\
40	133.150497\\
41	104.561329\\
42	127.287279\\
43	153.544953\\
44	125.139605\\
};
\addlegendentry{sym2IRKA}

\end{axis}
\end{tikzpicture}%
	\caption{Optimization times example 1}
	\label{pic:TimesEx1}
\end{figure}
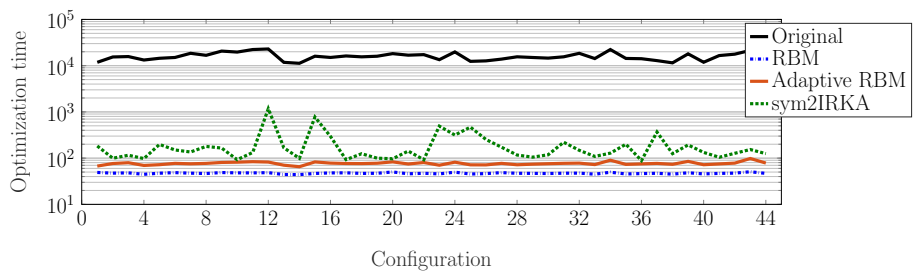

\subsection{Example 2}
The second example contains a mass oscillator with $2d+1 = n$ masses and $n+2$ springs.
We have two lines of $d$ consecutive masses $m_1,\dots,m_d$ and $m_{d+1},\dots,m_{2d}$ that are connected by springs.
The springs of the first line all have the stiffness value $k_1$ and the springs in the second line the stiffness value $k_2$, where the first masses are connected with these springs to a fixed object.
The last masses $m_d$ and $m_{2d}$ are connected to a mass $m_{2d+1}=m_n$ by springs with stiffness constants $k_1$ and $k_2$ while the mass $m_n$ is connected to a fixed object via a spring with a constant $k_1+k_2+k_3$.
This construction results in the stiffness matrix 
\begin{align*}
K = \begin{bmatrix}
K_{11} & & \kappa_1 \\
& K_{22} & \kappa_2 \\
\kappa_1^\T & \kappa_2^\T & k_1+k_2+k_3
\end{bmatrix},\quad K_{jj} = k_j \begin{bmatrix}
2  & -1     &        & & \\
-1 & 2      & -1     & &  \\
& \ddots & \ddots & \ddots  &\\
&        & -1     &2        & -1 \\
&        &        &-1 & 2
\end{bmatrix},\quad \kappa_j = \begin{bmatrix}
0\\
\vdots \\
0\\
k_j
\end{bmatrix},
\end{align*}

for $j=1, 2$.
We choose the dimension to be $d= 1000,~ n = 2001$ and set $k_1 = 400$, $k_2 =100$, $k_3 = 300$.
The $n=2d+1$ mass values are chosen as follows
\begin{align*}
m_j = \begin{cases}
100-\frac{j}{10},  & j = 1, \dots, 500,\\
\frac{j}{30} + 33,  & j = 501,\dots,1000,\\
100 - (j-99)\frac{5}{20} + \frac{(j-999)^2}{5000}, & j = 1001,\dots,2000,
\end{cases} \qquad m_{2001} = 100.
\end{align*}
The internal damping $D_{\intern}$ is built as described in Equation \eqref{eq:intDamp} with the scaling $\alpha=0.003$.
Additionally, there are disturbances, that effect $21$ masses. 
The effect to the masses is described by the matrix $B\in\R^{n\times 21}$ that consists of zero entries except the following
\begin{align*}
E(1:10,\, 1:10) &= \diag{1000,\, 900,\,\dots,\,100},\\
E(1001:1010,\, 11:20) &= \diag{1000,\, 900,\,\dots,\,100},\\
E(2001, 21) & = 2000.
\end{align*}
As output, we observe $42$ masses, or more detailed the displacements of the masses $490$ to $510$ and those of the masses on positions $1490$ to $1510$.
This is described by the output matrix $C\in\R^{42\times n}$:
\begin{align*}
C(490:510,\, 1:21) = I_{21},\qquad
C(1490:1510,\, 1:21) = I_{21}.
\end{align*}

In this example, we consider four damping values that are optimized.
We consider two dampers in the first row that are between the masses $m_j$ and $m_{j+5}$ and between the masses $m_{j+20}$ and $m_{j+25}$.
For the second row we follow the same pattern and add two dampers between the masses $m_k$ and $m_{k+5}$ and between $m_{k+20}$ and $m_{k+25}$.
Consequently, we optimize four damping values $g_1,~g_2,~g_3,~g_4$ that are saved in $g\in\R^4$.
The corresponding damping position matrix is then of the form
\[
F = \begin{bmatrix}
e_j - e_{j+5} & e_{j+20} -e_{j+25} & e_{k} - e_{k+5} & e_{k+20} - e_{k+25},
\end{bmatrix}
\]
where $j$ and $k$ are taken from the sets
\[
\left\{ (j, k)~ |~j \in \{250,\, 450,\, 650,\, 850\}, ~ k \in \{1150,\, 1250,\, 1350,\, 1450,\, 1550,\, 1650,\, 1750\}\right\}.
\]
This setting leads to $28$ different damping configurations.
We assume that the damping values $g_1,~g_2,~g_3,~g_4$  lie in the interval $[350, 7000]$.
For the optimization process, we set a tolerance of $5\cdot 10^{-4}$ and start at $g_0 = \begin{bmatrix}
1000 & 1000 & 1000 & 1000
\end{bmatrix}^{\T}$ for all damping configurations.
The tolerance for the function value error that indicates whether a bases is sufficiently detailed is set to be $\tol_f=10^{-2}$.
As test-parameter set $\cD_{\Test}$ for the RBM we use $21$ uniformly distributed parameters in $[350, 7000]^4$.
The first parameter $g_0^{\rr}$ that is used to obtain a first error equation basis is chosen to be $g_0^{\rr}=\begin{bmatrix}
100 & 100 & 100 & 100
\end{bmatrix}$ within the RBM and the adaptive RBM.

\begin{figure}[h!]
	\centering
%
%
\definecolor{mycolor1}{rgb}{1.0, 0.55, 0.0}%
\definecolor{mycolor2}{rgb}{0.0, 0.8, 0.6}%
\begin{tikzpicture}[scale=0.5]
\huge
\begin{axis}[%
width=4.521in,
height=3.566in,
at={(0.758in,0.603in)},
scale only axis,
xmin=0,
xmax=20,
xlabel style={font=\color{white!15!black},yshift=-0.8cm},
xlabel={test parameter},
xtick={0, 5,10,15,20},
ymode=linear,
ymin=0,
ymax=0.5e-01,
yminorticks=true,
ylabel style={font=\color{white!15!black},yshift=0.8cm},
ylabel={errors/error estimators},
ytick={0, 0.1e-01,0.2e-01,0.3e-01,0.4e-01, 0.5e-01},
axis background/.style={fill=white},
legend style={legend cell align=left, align=left, draw=white!15!black,xshift=8cm},
ymajorgrids,
yminorgrids,
xmajorgrids,
]
\addplot [color=mycolor1, dashdotted, line width=3.0pt]
  table[row sep=crcr]{%
1	0.000325361369224412\\
2	0.00355675718896567\\
3	0.00362667933559696\\
4	0.00316710240211202\\
5	0.00530589229141618\\
6	0.002721167279544\\
7	0.0038215087814118\\
8	0.00263997002884893\\
9	0.00563901219670457\\
10	0.00191586288829243\\
11	0.00665851532652345\\
12	0.0070240153302081\\
13	0.00521113195413332\\
14	0.00733740685876872\\
15	0.00319179603978845\\
16	0.00449645862150266\\
17	0.00395014223213252\\
18	0.00270628882250613\\
19	0.00764662584274312\\
20	0.00801259200818688\\
};
\addlegendentry{$\mathrm{tr}(\mathcal{CE}(g)\mathcal{C}^{\mathrm{T}})/\mathrm{tr}(\mathcal{C}\widetilde{P}(g)\mathcal{C}^{\mathrm{T}})$}

\addplot [color=blue, line width=3.0pt]
  table[row sep=crcr]{%
1	0.00215526624866213\\
2	0.00609555032909734\\
3	0.00583149299113374\\
4	0.00585627692075802\\
5	0.00801115209479173\\
6	0.00542664557801817\\
7	0.00649725149151078\\
8	0.00488354584932937\\
9	0.00783655168242209\\
10	0.0043981864671904\\
11	0.00875299444191907\\
12	0.00888978854200777\\
13	0.00735132435149531\\
14	0.00922821723505706\\
15	0.00520211017257472\\
16	0.00653722185330659\\
17	0.00590808071076411\\
18	0.00490726549997182\\
19	0.00922314846974634\\
20	0.00937924005482978\\
};
\addlegendentry{$\mathrm{tr}(\mathcal{C\widetilde{E}}(g)\mathcal{C}^{\mathrm{T}})/\mathrm{tr}(\mathcal{C}\widetilde{P}(g)\mathcal{C}^{\mathrm{T}})$}

\addplot [color=mycolor2, line width=3.0pt]
table[row sep=crcr]{%
	1	0.0127307073814356\\
	2	0.0255066121949011\\
	3	0.0205597963440585\\
	4	0.0319562553400312\\
	5	0.0357270314915193\\
	6	0.0333244396315991\\
	7	0.0359520493388318\\
	8	0.0292718362342471\\
	9	0.0371520233391344\\
	10	0.0376288549053617\\
	11	0.0495639770053865\\
	12	0.0443892404157463\\
	13	0.0488980151311791\\
	14	0.0493959629427177\\
	15	0.0317625410381581\\
	16	0.0382407050615745\\
	17	0.0412541180714153\\
	18	0.0441410408042541\\
	19	0.0534948068958103\\
	20	0.0487203798230225\\
};
\addlegendentry{$\|\widetilde{\mathcal{E}}_{11}(g)\|_{\mathrm{F}}/\|\widetilde{P}_{11}(g)\|_{\mathrm{F}}$}

\addplot [color=red, dashdotted, line width=3.0pt]
  table[row sep=crcr]{%
1	0.0124204601711904\\
2	0.02484630887451\\
3	0.0200364364383357\\
4	0.0311945535949302\\
5	0.0348269498088753\\
6	0.0324967491423689\\
7	0.0351053417284242\\
8	0.0286010674560917\\
9	0.0362476589469863\\
10	0.0368423512882056\\
11	0.0483600822762529\\
12	0.0433079340022074\\
13	0.0478388342789562\\
14	0.0482204088963016\\
15	0.0310452006417011\\
16	0.0373971273217252\\
17	0.0403400140860647\\
18	0.0432427908852498\\
19	0.052242854764713\\
20	0.0475785472936573\\
};
\addlegendentry{$\|\mathcal{E}_{11}(g)\|_{\mathrm{F}}/\|\widetilde{P}_{11}(g)\|_{\mathrm{F}}$}
\end{axis}

\end{tikzpicture}%
	\caption{Errors estimator for Example 2 for the first damping configuration and the fifth step of the RBM.}
	\label{pic:ErrorsEst2}
\end{figure}
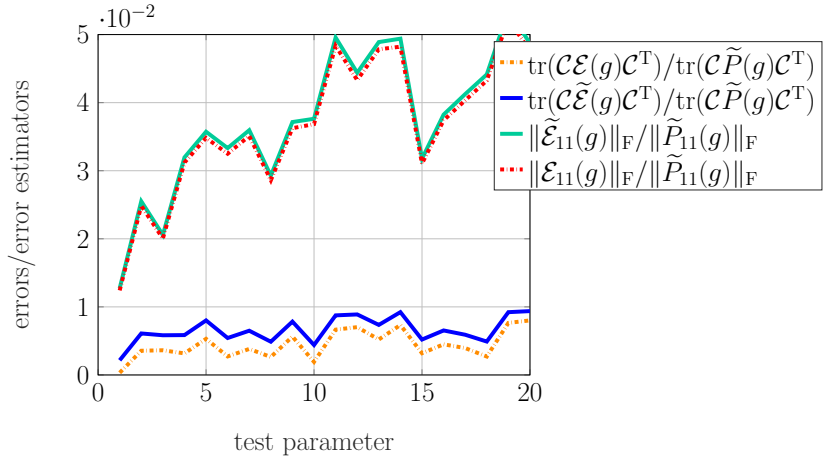

We evaluate the quality of the error estimate after the first step of the RBM for the fifth damper configuration in \Cref{pic:ErrorsEst2}.
We observe that the relative error in the position controllability Gramian and the corresponding estimate are very close, so the error is well approximated.
In this example, it can be observed that the error in the energy response is also well approximated.

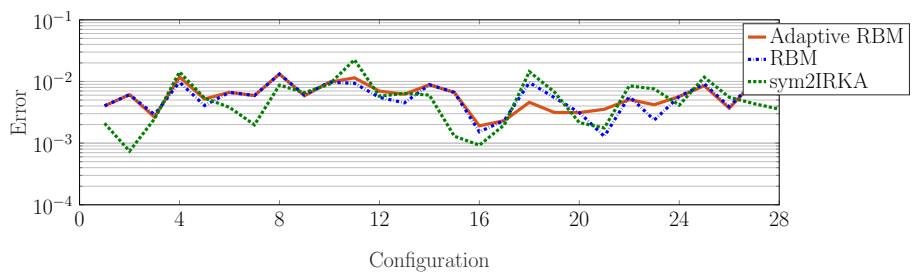
\begin{figure}[h!]
	\centering
%
%
\definecolor{mycolor1}{rgb}{0.85098,0.32549,0.09804}%
\definecolor{mycolor2}{rgb}{0.00000,0.49804,0.00000}%
\begin{tikzpicture}[scale=0.4]
\huge
\begin{axis}[%
width=9.155in,
height=2.422in,
at={(1.536in,0.584in)},
scale only axis,
xmin=0,
xmax=28,
xlabel style={font=\color{white!15!black},yshift=-0.8cm},
xlabel={Configuration},
xtick={0, 4,8,12,16,20,24,28},
ymode=log,
ymin=0.0001,
ymax=0.1,
yminorticks=true,
ylabel style={font=\color{white!15!black},yshift=0.8cm},
ylabel={Error},
axis background/.style={fill=white},
ymajorgrids,
yminorgrids,
legend style={legend cell align=left, align=left, draw=white!15!black,xshift=4.8cm}
]
\addplot [color=mycolor1, line width=3.0pt]
  table[row sep=crcr]{%
1	0.0040225332734859\\
2	0.00611031641500916\\
3	0.00263727852776901\\
4	0.0118316315273683\\
5	0.00519467247018509\\
6	0.00667167047885757\\
7	0.00591482969934829\\
8	0.0131918798010663\\
9	0.00586492977600293\\
10	0.00973515966789245\\
11	0.0114660882442601\\
12	0.00698226123735284\\
13	0.00621239587171511\\
14	0.00887206797281229\\
15	0.00663559252147376\\
16	0.00190764978079723\\
17	0.00231796914039479\\
18	0.00459587744303437\\
19	0.00314484688806571\\
20	0.00312871013600977\\
21	0.00353942881032631\\
22	0.00508842302181157\\
23	0.00419380292853511\\
24	0.00571206616523281\\
25	0.00850763486539479\\
26	0.00368122116487034\\
27	0.00995929072812766\\
28	0.00721795851045785\\
};
\addlegendentry{Adaptive RBM}

\addplot [color=blue, dashdotted, line width=3.0pt]
  table[row sep=crcr]{%
1	0.0040225332734859\\
2	0.00611031641500916\\
3	0.0029118136046656\\
4	0.00977504720390863\\
5	0.00407286520483841\\
6	0.00667167047885757\\
7	0.00591482969934829\\
8	0.0131918798010663\\
9	0.00606530853215923\\
10	0.00973515966789245\\
11	0.00935030051036397\\
12	0.00552766924833249\\
13	0.00449984721613377\\
14	0.00887206797281229\\
15	0.00663559252147376\\
16	0.00152619790217891\\
17	0.00233677691259487\\
18	0.00961341272845086\\
19	0.00548382393685794\\
20	0.00312871013600977\\
21	0.00128828186960954\\
22	0.0057515193629914\\
23	0.00236728293073784\\
24	0.00571206616523281\\
25	0.0091544859671864\\
26	0.0039030365548755\\
27	0.00995929072812766\\
28	0.00721795851045785\\
};
\addlegendentry{RBM}

\addplot [color=mycolor2, dotted, line width=3.0pt]
  table[row sep=crcr]{%
1	0.00209916040202903\\
2	0.000741889661440154\\
3	0.00247756098431474\\
4	0.0142526873790428\\
5	0.00537689023903715\\
6	0.0038013315424855\\
7	0.00198277277997369\\
8	0.00872742060601473\\
9	0.00658043597600811\\
10	0.00905172251662816\\
11	0.022480399073334\\
12	0.00578528823603053\\
13	0.00640356516587544\\
14	0.00602173421796108\\
15	0.00129542251592613\\
16	0.000927122973361299\\
17	0.00197877363751826\\
18	0.014409399353678\\
19	0.00676896882645772\\
20	0.00215102333770369\\
21	0.0017634223769724\\
22	0.00849783536964515\\
23	0.00756725552953131\\
24	0.0041335541377362\\
25	0.0117384005491257\\
26	0.0055403478691654\\
27	0.00436182839442496\\
28	0.00355474372316054\\
};
\addlegendentry{sym2IRKA}

\end{axis}

\begin{axis}[%
width=5.833in,
height=4.375in,
at={(0in,0in)},
scale only axis,
xmin=0,
xmax=1,
ymin=0,
ymax=1,
axis line style={draw=none},
ticks=none,
axis x line*=bottom,
axis y line*=left,
legend style={legend cell align=left, align=left, draw=white!15!black}
]
\end{axis}
\end{tikzpicture}%
	\caption{Errors example 2}
	\label{pic:ErrorsEx2}
\end{figure}

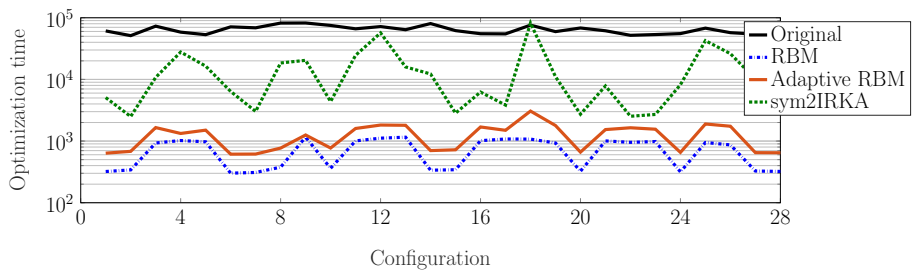
\begin{figure}[h!]
	\centering
%
%
\definecolor{mycolor1}{rgb}{0.85098,0.32549,0.09804}%
\definecolor{mycolor2}{rgb}{0.00000,0.49804,0.00000}%
\begin{tikzpicture}[scale=0.4]
\huge
\begin{axis}[%
width=9.155in,
height=2.422in,
at={(1.536in,0.584in)},
scale only axis,
xmin=0,
xmax=28,
xlabel style={font=\color{white!15!black},yshift=-0.8cm},
xlabel={Configuration},
xtick={0, 4,8,12,16,20,24,28},
ymode=log,
ymin=100,
ymax=100000,
yminorticks=true,
ylabel style={font=\color{white!15!black},yshift=0.8cm},
ylabel={Optimization time},
axis background/.style={fill=white},
ymajorgrids,
yminorgrids,
legend style={legend cell align=left, align=left, draw=white!15!black,xshift=4.8cm}
]
\addplot [color=black, line width=3.0pt]
  table[row sep=crcr]{%
1	61175.028364\\
2	51292.083801\\
3	73110.708303\\
4	58448.737637\\
5	53186.198345\\
6	71202.42686\\
7	68715.113676\\
8	81976.120628\\
9	82271.572102\\
10	75091.128436\\
11	66189.197914\\
12	71962.21188\\
13	63997.934492\\
14	80528.330508\\
15	61865.271622\\
16	55096.786346\\
17	54762.661047\\
18	75942.337884\\
19	59592.222232\\
20	68264.91785\\
21	61509.163274\\
22	51916.607528\\
23	53424.477474\\
24	55334.84365\\
25	67838.356153\\
26	57303.082269\\
27	54257.790897\\
28	75297.017262\\
};
\addlegendentry{Original}

\addplot [color=blue, dashdotted, line width=3.0pt]
  table[row sep=crcr]{%
1	320.79036\\
2	340.704647\\
3	930.552133\\
4	1020.12831\\
5	972.533083\\
6	301.549975\\
7	309.751066\\
8	377.190795\\
9	1159.458775\\
10	361.741871\\
11	1006.977068\\
12	1114.300104\\
13	1156.137044\\
14	336.71553\\
15	343.068004\\
16	1014.893946\\
17	1080.641099\\
18	1075.958242\\
19	931.0514\\
20	324.290509\\
21	1005.315448\\
22	953.535362\\
23	981.820854\\
24	316.056816\\
25	947.484993\\
26	866.340082\\
27	326.760746\\
28	321.624548\\
};
\addlegendentry{RBM}

\addplot [color=mycolor1, line width=3.0pt]
  table[row sep=crcr]{%
1	636.247601\\
2	681.031067\\
3	1654.801518\\
4	1328.694932\\
5	1503.222975\\
6	613.119386\\
7	615.463395\\
8	769.680041\\
9	1254.770102\\
10	771.024157\\
11	1600.575699\\
12	1819.513545\\
13	1803.14917\\
14	701.070457\\
15	722.008268\\
16	1692.890845\\
17	1494.780696\\
18	3064.036461\\
19	1787.704555\\
20	661.850879\\
21	1536.202709\\
22	1643.155974\\
23	1557.282912\\
24	655.332571\\
25	1892.285588\\
26	1740.104262\\
27	648.551979\\
28	646.084758\\
};
\addlegendentry{Adaptive RBM}

\addplot [color=mycolor2, dotted, line width=3.0pt]
  table[row sep=crcr]{%
1	5074.722803\\
2	2499.658954\\
3	10778.479585\\
4	27786.107609\\
5	16445.712692\\
6	6391.064815\\
7	3026.798811\\
8	18485.261926\\
9	20395.458609\\
10	4381.192847\\
11	24620.119443\\
12	56794.237446\\
13	15942.697455\\
14	12257.421989\\
15	2837.364587\\
16	6236.366867\\
17	3831.677968\\
18	82091.418795\\
19	11084.358709\\
20	2726.800569\\
21	7827.220416\\
22	2525.913843\\
23	2703.304284\\
24	8209.584637\\
25	42568.306881\\
26	26152.159513\\
27	10107.189021\\
28	12480.936928\\
};
\addlegendentry{sym2IRKA}

\end{axis}
\end{tikzpicture}%
	\caption{Optimization times example 2}
	\label{pic:TimesEx2}
\end{figure}

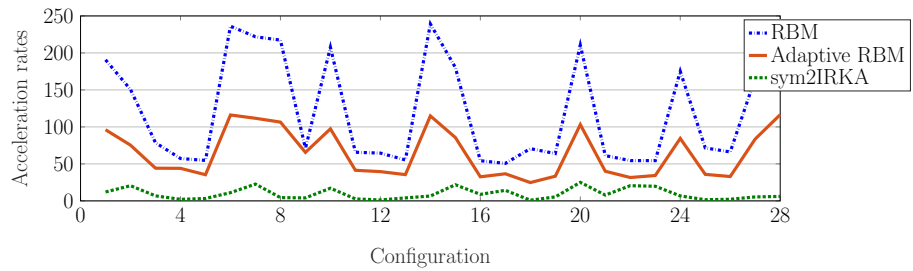
\begin{figure}[h!]
	\centering
%
%
\definecolor{mycolor1}{rgb}{0.85098,0.32549,0.09804}%
\definecolor{mycolor2}{rgb}{0.00000,0.49804,0.00000}%
\begin{tikzpicture}[scale=0.4]
\huge

\begin{axis}[%
width=9.155in,
height=2.422in,
at={(1.536in,0.584in)},
scale only axis,
xmin=0,
xmax=28,
xlabel style={font=\color{white!15!black},yshift=-0.8cm},
xlabel={Configuration},
xtick={0, 4,8,12,16,20,24,28},
ymin=0,
ymax=250,
ylabel style={font=\color{white!15!black},yshift=0.8cm},
ylabel={Acceleration rates},
axis background/.style={fill=white},
ymajorgrids,
legend style={legend cell align=left, align=left, draw=white!15!black,xshift=4.8cm}
]
\addplot [color=blue, dashdotted, line width=3.0pt]
  table[row sep=crcr]{%
1	190.700956113519\\
2	150.547062544175\\
3	78.5670202778419\\
4	57.2954765239286\\
5	54.6883178317544\\
6	236.12148155542\\
7	221.839797238987\\
8	217.333301116216\\
9	70.9568756353584\\
10	207.582075661902\\
11	65.7305911101443\\
12	64.5806382155736\\
13	55.3549726860927\\
14	239.158349803468\\
15	180.329470835759\\
16	54.2882205211223\\
17	50.6760857954376\\
18	70.5811200840265\\
19	64.0052979158831\\
20	210.505444826324\\
21	61.1839432054564\\
22	54.4464417335369\\
23	54.4136715535684\\
24	175.078785992706\\
25	71.5983436721303\\
26	66.1438659708694\\
27	166.0474569274\\
28	234.114646192989\\
};
\addlegendentry{RBM}

\addplot [color=mycolor1, line width=3.0pt]
  table[row sep=crcr]{%
1	96.1497194926162\\
2	75.3153362400133\\
3	44.1809531280597\\
4	43.9895842373846\\
5	35.3814432253472\\
6	116.131423154837\\
7	111.647766925277\\
8	106.506751196891\\
9	65.5670484743507\\
10	97.3914082383284\\
11	41.3533692629179\\
12	39.5502479647658\\
13	35.4923128694893\\
14	114.864818084896\\
15	85.6849905519365\\
16	32.5459769061484\\
17	36.6359166890124\\
18	24.7850633798316\\
19	33.3344914657892\\
20	103.142444946424\\
21	40.0397440478671\\
22	31.5956661141653\\
23	34.3062118400744\\
24	84.4378047096945\\
25	35.8499565727285\\
26	32.93083266352\\
27	83.6598956658183\\
28	116.543559230661\\
};
\addlegendentry{Adaptive RBM}

\addplot [color=mycolor2, dotted, line width=3.0pt]
  table[row sep=crcr]{%
1	12.054851218245\\
2	20.5196327758719\\
3	6.78302609625437\\
4	2.10352376300696\\
5	3.23404642541715\\
6	11.1409333062757\\
7	22.7022402104413\\
8	4.43467455079435\\
9	4.0338181984148\\
10	17.1394255076944\\
11	2.68841904147703\\
12	1.26706889846741\\
13	4.01424756837048\\
14	6.5697608012735\\
15	21.8037794316068\\
16	8.83475708229209\\
17	14.2920833912314\\
18	0.92509471755683\\
19	5.376244471736\\
20	25.0348040212691\\
21	7.85836606163104\\
22	20.5535939683276\\
23	19.7626577926142\\
24	6.74027324118323\\
25	1.59363529168878\\
26	2.1911415093853\\
27	5.3682374777267\\
28	6.03296192396238\\
};
\addlegendentry{sym2IRKA}

\end{axis}

\end{tikzpicture}%
	\caption{Acceleration rates example 2}
	\label{pic:AccelerationsEx2}
\end{figure}

We compare again the results generated by the methods presented above with the symmetric IRKA approach for second order systems (\texttt{sym2IRKA}) \cite{morTomBG18}.
In every iteration, $120$ vectors are added to the basis and the method stops if the relative error between two consecutive optima is smaller that the tolerance $\tol_{\mathrm{IRKA}}=10^{-3}$ .
The relative errors between the optimal damping gain and the approximations we obtain using the methods presented in the previous sections are shown in \Cref{pic:ErrorsEx2}.
Additionally, we evaluate the optimization times.
Outside of the applied methods, we determine a low-rank factor of the solution of the Lyapunov equation \eqref{eq:Lyap} for the undamped system. This low-rank factor is included into all the bases and is not taken into account in the time measures. 
This solving takes $54$ seconds.
We compare the optimization times and the acceleration rates for the different methods in \Cref{pic:TimesEx2} and 
\Cref{pic:AccelerationsEx2}, respectively.
They are in average $89$ for the reduced basis method and $51$ for the adaptive procedure.
These results are comparable with those from \cite{morTomBG18}.
Here the acceleration rate is in average $4$ in our implementation, that is less improvement than with our method.
However, in \cite{morTomBG18}, acceleration rates of up to $208$ were achieved for this example, which we could not reproduce.

\begin{figure}[h!]
	\centering
%
%
\definecolor{mycolor1}{rgb}{0.00000,0.44700,0.74100}%
\begin{tikzpicture}[scale=0.4]
\huge

\begin{axis}[%
width=9.155in,
height=2.422in,
at={(1.536in,0.584in)},
scale only axis,
xmin=0,
xmax=30,
xlabel style={font=\color{white!15!black},yshift=-0.8cm},
xlabel={Configuration},
xtick={0, 4,8,12,16,20,24,28,32,36,40,44},
ymin=750,
ymax=1050,
ylabel style={font=\color{white!15!black},yshift=0.8cm},
ylabel={Function value},
axis background/.style={fill=white},
ymajorgrids
]
\addplot [color=mycolor1, line width=3.0pt, forget plot]
  table[row sep=crcr]{%
1	794.262382871874\\
2	782.891558348899\\
3	859.920981020981\\
4	763.515344941595\\
5	825.390961040074\\
6	894.605544319072\\
7	882.499252987989\\
8	897.828552577089\\
9	894.865743632487\\
10	935.774143799376\\
11	865.461735681554\\
12	903.325996462364\\
13	974.320454126973\\
14	1001.99741802361\\
15	837.023361509735\\
16	829.090694002077\\
17	899.050182294521\\
18	805.73748950947\\
19	859.999621960504\\
20	921.599991570255\\
21	932.310557516001\\
22	785.172111631647\\
23	786.401929360219\\
24	857.324055712903\\
25	761.695155605984\\
26	807.809210772613\\
27	887.282585450035\\
28	874.577912901851\\
};
\end{axis}
\end{tikzpicture}%
	\caption{Function values obtained in Example 2}
	\label{pic:FunValEx2}
\end{figure}
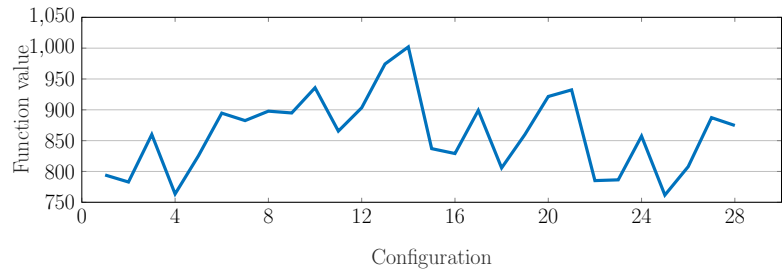

In \Cref{pic:FunValEx2} the function values for all $28$ damping configurations are evaluated. We see that the optimal damping configuration is the $25$-th one that has the damping positions $j=850,~k=1450$.

\section{Conclusion}\label{sec:Conclusion}
In this article, we have dealt with the damping optimization problem in a large-scale vibrational system.
We have presented two different approaches based on the reduced basis method (RBM), which can be used to find optimal damping values in a reasonable time.
In the first method, the RBM for Lyapunov equations has been integrated into the calculation of the energy response, and, subsequently a reduced objective functional has been optimized.
In the second method, we have integrated the RBM into the optimization process so that the basis is only enriched when necessary.
In addition, we derived new error estimators that were used to evaluate the quality of the reduced energy response and the reduced position controllability Gramian.
Finally, we have applied our methods to some mechanical systems and have observed that our approaches improve on or at are least comparable to existing methods.

\section*{Declarations}
\textbf{Conflict of interest} The authors declare no competing interests.

\bibliographystyle{plainurl}
\bibliography{csc,mor,software,References}

\begin{thebibliography}{10}

\bibitem{BarS72}
R.~H. Bartels and G.~W. Stewart.
\newblock Solution of the matrix equation ${AX}+{XB}={C}$: {A}lgorithm 432.
\newblock {\em Comm. {ACM}}, 15:820--826, 1972.
\newblock \href {https://doi.org/10.1145/361573.361582}
  {\path{doi:10.1145/361573.361582}}.

\bibitem{morBeaGT20}
C.~Beattie, S.~Gugercin, and Z.~Tomljanovi{\'c}.
\newblock Sampling-free model reduction of systems with low-rank
  parameterization.
\newblock {\em Adv. Comput. Math.}, 46(6):83, 2020.
\newblock \href {https://doi.org/10.1007/s10444-020-09825-8}
  {\path{doi:10.1007/s10444-020-09825-8}}.

\bibitem{morBenKTetal16}
P.~Benner, P.~K{\"u}rschner, Z.~Tomljanovi{\'c}, and N.~Truhar.
\newblock Semi-active damping optimization of vibrational systems using the
  parametric dominant pole algorithm.
\newblock {\em Z. Angew. Math. Mech.}, 96(5):604--619, 2016.
\newblock \href {https://doi.org/10.1002/zamm.201400158}
  {\path{doi:10.1002/zamm.201400158}}.

\bibitem{BenQ99}
P.~Benner and E.~S. Quintana-Ort{\'\i}.
\newblock Solving stable generalized {L}yapunov equations with the matrix sign
  function.
\newblock {\em Numer. Algorithms}, 20(1):75--100, 1999.
\newblock \href {https://doi.org/10.1023/A:1019191431273}
  {\path{doi:10.1023/A:1019191431273}}.

\bibitem{morBenTT11}
P.~Benner, Z.~Tomljanovi{\'c}, and N.~Truhar.
\newblock Damping optimization for linear vibrating systems using dimension
  reduction.
\newblock In J.~N{\'a}prstek, J.~Hor{\'a}{\v c}ek, M.~Okrouhl{\'\i}k,
  B.~Marvalov{\'a}, F.~Verhulst, and J.~T. Sawicki, editors, {\em Vibration
  Problems ICOVP 2011}, volume 139, Part 5 of {\em Springer Proceedings in
  Physics}, pages 297--305, Prag, Czech Republic, 2011. Springer-Verlag.
\newblock \href {https://doi.org/10.1007/978-94-007-2069-5_41}
  {\path{doi:10.1007/978-94-007-2069-5_41}}.

\bibitem{morBenTT11a}
P.~Benner, Z.~Tomljanovi{\'c}, and N.~Truhar.
\newblock Dimension reduction for damping optimization in linear vibrating
  systems.
\newblock {\em Z. Angew. Math. Mech.}, 91(3):179--191, 2011.
\newblock \href {https://doi.org/10.1002/zamm.201000077}
  {\path{doi:10.1002/zamm.201000077}}.

\bibitem{morBenTT13}
P.~Benner, Z.~Tomljanovi{\'c}, and N.~Truhar.
\newblock Optimal damping of selected eigenfrequencies using dimension
  reduction.
\newblock {\em Numer. Lin. Alg. Appl.}, 20(1):1--17, 2013.
\newblock \href {https://doi.org/10.1002/nla.833} {\path{doi:10.1002/nla.833}}.

\bibitem{morCheFB21}
S.~Chellappa, L.~Feng, and P.~Benner.
\newblock A training set subsampling strategy for the reduced basis method.
\newblock {\em J. Sci. Comput.}, 89(63), 2021.
\newblock \href {https://doi.org/10.1007/s10915-021-01665-y}
  {\path{doi:10.1007/s10915-021-01665-y}}.

\bibitem{morCheFdetal21}
S.~Chellappa, L.~Feng, V.~de~la Rubia, and P.~Benner.
\newblock Inf-sup-constant-free state error estimator for model order reduction
  of parametric systems in electromagnetics.
\newblock e-prints 2104.12802, arXiv, 2021.
\newblock math.NA (accepted in {IEEE} Trans. Microw. Theory Techn.).
\newblock URL: \url{https://arxiv.org/abs/2104.12802}.

\bibitem{morDen19}
J.~Deni{\ss}en.
\newblock {\em On Vibrational Analysis and Reduction for Damped Linear
  Systems}.
\newblock {Dissertation}, Department of Mathematics, Otto-von-Guericke
  University, Magdeburg, Germany, 2019.
\newblock URL: \url{https://opendata.uni-halle.de//handle/1981185920/14371}.

\bibitem{DymD21}
A.~Dymarek and T.~Dzitkowski.
\newblock The use of synthesis methods in position optimisation and selection
  of tuned mass damper (tmd) parameters for systems with many degrees of
  freedom.
\newblock {\em Archives of Control Sciences}, 31(LXVII)(1):185--21, 2021.
\newblock \href {https://doi.org/10.24425/acs.2021.136886}
  {\path{doi:10.24425/acs.2021.136886}}.

\bibitem{morEftKP11}
J.~L. Eftang, D.~J. Knezevic, and A.~T. Patera.
\newblock An hp certified reduced basis method for parametrized parabolic
  partial differential equations.
\newblock {\em Math. Comput. Model. Dyn. Syst.}, 17(4):395--422, 2011.
\newblock \href {https://doi.org/10.1080/13873954.2011.547670}
  {\path{doi:10.1080/13873954.2011.547670}}.

\bibitem{morFenB21}
L.~Feng and P.~Benner.
\newblock On error estimation for reduced-order modeling of linear
  non-parametric and parametric systems.
\newblock {\em {ESAIM}: Math. Model. Numer. Anal.}, 55(2):561--594, 2021.
\newblock \href {https://doi.org/10.1051/m2an/2021001}
  {\path{doi:10.1051/m2an/2021001}}.

\bibitem{Gen09}
G.~Genta.
\newblock {\em Vibration Dynamics and Control}.
\newblock Mechanical Engineering Series. Springer, 2009.

\bibitem{morGugAB08}
S.~Gugercin, A.~C. Antoulas, and C.~Beattie.
\newblock {$\mathcal{H}_2$} model reduction for large-scale linear dynamical
  systems.
\newblock {\em {SIAM} J. Matrix Anal. Appl.}, 30(2):609--638, 2008.
\newblock \href {https://doi.org/10.1137/060666123}
  {\path{doi:10.1137/060666123}}.

\bibitem{GurM92}
M.~G\"{u}rg\"{o}ze and P.C. M\"{u}ller.
\newblock Optimal positioning of dampers in multi-body systems.
\newblock {\em Journal of Sound and Vibration}, 158:517--530, 1992.
\newblock \href {https://doi.org/10.1016/0022-460X(92)90422-T}
  {\path{doi:10.1016/0022-460X(92)90422-T}}.

\bibitem{Ham82b}
S.~J. Hammarling.
\newblock Numerical solution of the stable, non-negative definite {L}yapunov
  equation.
\newblock {\em {IMA} J. Numer. Anal.}, 2:303--323, 1982.

\bibitem{morHesRS16}
J.~S. Hesthaven, G.~Rozza, and B.~Stamm.
\newblock {\em Certified {R}educed {B}asis {M}ethods for {P}arametrized
  {P}artial {D}ifferential {E}quations}.
\newblock SpringerBriefs in Mathematics. Springer, Cham, 2016.
\newblock \href {https://doi.org/10.1007/978-3-319-22470-1}
  {\path{doi:10.1007/978-3-319-22470-1}}.

\bibitem{morHesSZ14}
J.~S. Hesthaven, B.~Stamm, and S.~Zhang.
\newblock Efficient greedy algorithms for high-dimensional parameter spaces
  with applications to empirical interpolation and reduced basis methods.
\newblock {\em ESAIM: Math. Model. Numer. Anal.}, 48(1):259--283, 2014.
\newblock \href {https://doi.org/10.1051/m2an/2013100}
  {\path{doi:10.1051/m2an/2013100}}.

\bibitem{Inm06}
D.~J. Inman.
\newblock {\em Vibration with Control}.
\newblock John Wiley \& Sons Ltd., Virginia Tech, USA, 2006.

\bibitem{JakMTetal21}
N.~Jakov{\v c}evi{\'c}~Stor, T.~Mitchell, Z.~Tomljanovi{\'c}, and M.~Ugrica.
\newblock Fast optimization of viscosities for frequency-weighted damping of
  second-order systems.
\newblock e-print arxiv:2104.04035, arXiv, April 2021.
\newblock math.NA.
\newblock URL: \url{https://arxiv.org/abs/2104.04035}.

\bibitem{Kan13}
Y.~Kanno.
\newblock {D}amper placement optimization in a shear building model
  withdiscrete design variables: a mixed-integer second-order coneprogramming
  approach.
\newblock {\em Earthquake Engng Struct. Dyn}, 42:1657--1676, 2013.
\newblock \href {https://doi.org/10.1002/eqe.2292}
  {\path{doi:10.1002/eqe.2292}}.

\bibitem{KanPTetal19}
Y.~Kanno, M.~Puva{\v c}a, Z.~Tomljanovi{\'c}, and N.~Truhar.
\newblock Optimization of damping positions in a mechanical system.
\newblock {\em Rad Hrvat. Akad. Znan. Umjet. Mat. Znan.}, 23:141--157, 2019.
\newblock \href {https://doi.org/10.21857/y26kec33q9}
  {\path{doi:10.21857/y26kec33q9}}.

\bibitem{Kue16}
Patrick K{\"u}rschner.
\newblock {\em Efficient Low-Rank Solution of Large-Scale Matrix Equations}.
\newblock {D}issertation, Otto-von-Guericke-Universit{\"a}t, Magdeburg,
  Germany, April 2016.
\newblock URL: \url{http://hdl.handle.net/11858/00-001M-0000-0029-CE18-2}.

\bibitem{KuzTT12}
I.~Kuzmanovi{\'c}, Z.~Tomljanovi{\'c}, and N.~Truhar.
\newblock Optimization of material with modal damping.
\newblock {\em Appl. Math. Comput.}, 218(13):7326--7338, 2012.
\newblock \href {https://doi.org/10.1016/j.amc.2012.01.011}
  {\path{doi:10.1016/j.amc.2012.01.011}}.

\bibitem{MueS85}
P.~C. M{\"u}ller and W.~O. Schiehlen.
\newblock {\em {L}inear vibrations: {A} theoretical treatment of
  multi-degree-of-freedom vibrating systems}.
\newblock Martinus Hijhoff publishers, 1985.

\bibitem{PazK18}
M.~Paz and Y.~H. Kim.
\newblock {\em {S}tructural {D}ynamics: {T}heory and {C}omputation}.
\newblock Springer, 2018.

\bibitem{Pen00b}
T.~Penzl.
\newblock A cyclic low rank {S}mith method for large sparse {L}yapunov
  equations.
\newblock {\em {SIAM} J. Sci. Comput.}, 21(4):1401--1418, 2000.
\newblock \href {https://doi.org/10.1137/S1064827598347666}
  {\path{doi:10.1137/S1064827598347666}}.

\bibitem{morPrzV21}
J.~Przybilla and M.~Voigt.
\newblock Model reduction of parametric differential-algebraic systems by
  balanced truncation.
\newblock e-print 2108.08646, arXiv, 2021.
\newblock math.DS.
\newblock URL: \url{https://arxiv.org/abs/2108.08646}.

\bibitem{morQuaMN16}
A.~Quarteroni, A.~Manzoni, and F.~Negri.
\newblock {\em Reduced Basis Methods for Partial Differential Equations},
  volume~92 of {\em La Matematica per il 3+2}.
\newblock Springer International Publishing, 2016.
\newblock ISBN: 978-3-319-15430-5.

\bibitem{morRob80}
J.~D. Roberts.
\newblock Linear model reduction and solution of the algebraic {R}iccati
  equation by use of the sign function.
\newblock {\em Internat. J. Control}, 32(4):677--687, 1980.
\newblock (Reprint of Technical Report No. TR-13, CUED/B-Control, Cambridge
  University, Engineering Department, 1971).
\newblock \href {https://doi.org/10.1080/00207178008922881}
  {\path{doi:10.1080/00207178008922881}}.

\bibitem{morRozHP08}
G.~Rozza, D.~B.~P. Huynh, and A.~T. Patera.
\newblock Reduced basis approximation and a posteriori error estimation for
  affinely parametrized elliptic coercive partial differential equations.
\newblock {\em Archives of Computational Methods in Engineering},
  15(3):229--275, 2008.
\newblock \href {https://doi.org/10.1007/s11831-008-9019-9}
  {\path{doi:10.1007/s11831-008-9019-9}}.

\bibitem{morSchH15}
A.~Schmidt and B.~Haasdonk.
\newblock Reduced basis approximation of large scale algebraic {R}iccati
  equations.
\newblock Technical report, University of Stuttgart, 2015.
\newblock URL:
  \url{http://www.simtech.uni-stuttgart.de/publikationen/prints.php?ID=999}.

\bibitem{morSchH18}
A.~Schmidt and B.~Haasdonk.
\newblock Reduced basis approximation of large scale algebraic {R}iccati
  equations.
\newblock {\em ESAIM: Control Optim. Calculus Variations}, 24(1):129 -- 151,
  2018.
\newblock \href {https://doi.org/10.1051/cocv/2017011}
  {\path{doi:10.1051/cocv/2017011}}.

\bibitem{SimD09}
V.~Simoncini and V.~Druskin.
\newblock Convergence analysis of projection methods for the numerical solution
  of large {L}yapunov equations.
\newblock {\em {SIAM} J. Numer. Anal.}, 47(2):828--843, 2009.
\newblock \href {https://doi.org/10.1137/070699378}
  {\path{doi:10.1137/070699378}}.

\bibitem{morSonS17}
N.~T. Son and T.~Stykel.
\newblock Solving parameter-dependent {L}yapunov equations using the reduced
  basis method with application to parametric model order reduction.
\newblock {\em {SIAM} J. Matrix Anal. Appl.}, 38(2):478--504, 2017.
\newblock \href {https://doi.org/10.1137/15M1027097}
  {\path{doi:10.1137/15M1027097}}.

\bibitem{Tak97}
I.~Takewaki.
\newblock {O}ptimal damper placement for minimum transfer functions.
\newblock {\em Earthquake Engng Struct. Dyn.}, 26:1113--1997, 1997.
\newblock \href
  {https://doi.org/10.1002/(SICI)1096-9845(199711)26:11<1113::AID-EQE696>3.0.CO;2-X}
  {\path{doi:10.1002/(SICI)1096-9845(199711)26:11<1113::AID-EQE696>3.0.CO;2-X}}.

\bibitem{morTom11}
Z.~Tomljanovi{\'c}.
\newblock {\em Optimal damping for vibrating systems using dimension
  reduction}.
\newblock PhD thesis, Department of Mathematics, University of Osijek, Osijek,
  Croatia, 2011.
\newblock URL: \url{https://hrcak.srce.hr/68645}.

\bibitem{morTomBG18}
Z.~Tomljanovi{\'c}, C.~Beattie, and S.~Gugercin.
\newblock Damping optimization of parameter dependent mechanical systems by
  rational interpolation.
\newblock {\em Adv. Comput. Math.}, 44(6):1797--1820, 2018.
\newblock \href {https://doi.org/10.1007/s10444-018-9605-9}
  {\path{doi:10.1007/s10444-018-9605-9}}.

\bibitem{morTomV20}
Z.~Tomljanovi{\'c} and M.~Voigt.
\newblock Semi-active $\mathcal{H}_\infty$ damping optimization by adaptive
  interpolation.
\newblock {\em Numer. Lin. Alg. Appl.}, 27(4):e2300, 2020.
\newblock \href {https://doi.org/10.1002/nla.2300}
  {\path{doi:10.1002/nla.2300}}.

\bibitem{morTru04}
N.~Truhar.
\newblock An efficient algorithm for damper optimization for linear vibrating
  systems using {L}yapunov equation.
\newblock {\em J. Comput. Appl. Math.}, 127:169--182, 2004.
\newblock \href {https://doi.org/10.1016/j.cam.2004.02.005}
  {\path{doi:10.1016/j.cam.2004.02.005}}.

\bibitem{TruV07}
N.~Truhar and K.~Veseli{\'c}.
\newblock Bounds on the trace of a solution to the {L}yapunov equation with a
  general stable matrix.
\newblock {\em Systems Control Lett.}, 56(7--8):493--503, 2007.
\newblock \href {https://doi.org/10.1016/j.sysconle.2007.02.003}
  {\path{doi:10.1016/j.sysconle.2007.02.003}}.

\bibitem{morTruV09}
N.~Truhar and K.~Veseli{\'c}.
\newblock An efficient method for estimating the optimal dampers' viscosity for
  linear vibrating systems using {L}yapunov equation.
\newblock {\em {SIAM} J. Matrix Anal. Appl.}, 31(1):18--39, 2009.
\newblock \href {https://doi.org/10.1137/070683052}
  {\path{doi:10.1137/070683052}}.

\bibitem{morVerP05}
K.~Veroy and A.~T. Patera.
\newblock Certified real-time solution of the parametrized steady
  incompressible {N}avier-{S}tokes equations: rigorous reduced-basis a
  posteriori error bounds.
\newblock {\em Internat. J. Numer. Methods Fluids}, 47(8-9):773--788, 2005.
\newblock \href {https://doi.org/10.1002/fld.867} {\path{doi:10.1002/fld.867}}.

\bibitem{morVerPP03}
K.. Veroy, C.~Prud'homme, and A.~T. Patera.
\newblock Reduced-basis approximation of the viscous burgers equation: rigorous
  a posteriori error bounds.
\newblock {\em Comptes Rendus Mathematique. Academie des Sciences. Paris},
  337(9):619--624, 2003.
\newblock \href {https://doi.org/10.1016/j.crma.2003.09.023}
  {\path{doi:10.1016/j.crma.2003.09.023}}.

\bibitem{morVerPRetal03}
K.~Veroy, C.~Prud'Homme, D.~V. Rovas, and A.~T. Patera.
\newblock A posteriori error bounds for reduced-basis approximation of
  parametrized noncoercive and nonlinear elliptic partial differential
  equations.
\newblock In {\em 16th AIAA Computational Fluid Dynamics Conference}, Orlando,
  United States, 2003.
\newblock URL: \url{https://hal.archives-ouvertes.fr/hal-01219051}.

\bibitem{Ves11}
K.~Veseli{\'c}.
\newblock {\em Damped {O}scillations of {L}inear {S}ystems}, volume 2023 of
  {\em Lecture Notes in Math.}
\newblock Springer-Verlag, 2011.
\newblock \href {https://doi.org/10.1007/978-3-642-21335-9}
  {\path{doi:10.1007/978-3-642-21335-9}}.

\bibitem{ZhoDG96}
K.~Zhou, J.~C. Doyle, and K.~Glover.
\newblock {\em Robust and Optimal Control}.
\newblock Prentice-Hall, Upper Saddle River, NJ, 1996.
\newblock \href {https://doi.org/10.1007/978-1-4471-6257-5}
  {\path{doi:10.1007/978-1-4471-6257-5}}.

\end{thebibliography}

\end{document}